\documentclass[10pt]{article}
\usepackage{amsmath,amssymb,amsthm,amscd}
\numberwithin{equation}{section}

\def\p{\partial}
\def\b{\bar}

\def\L{\Lambda}

\def\cD{{\cal D}}

\def\cH{{\cal H}}

\def\M{{\mathfrak M}}

\def\cD{\mathcal D}

\def\cF{{\mathcal F}}

\def\cH{{\mathcal H}}

\def\cR{{\mathcal R}}
\def\cS{{\mathcal S}}

\def\cU{{\mathcal U}}
\def\cV{{\mathcal V}}

\def\cX{{\mathcal X}}

\def\bE{{\mathbf E}}

\def\RR{{\mathbb R}}
\def\CC{{\mathbb C}}
\def\NN{{\mathbb N}}

\newtheorem{prop}{Proposition}[section]
\newtheorem{theo}[prop]{Theorem}
\newtheorem{lem}[prop]{Lemma}
\newtheorem{cor}[prop]{Corollary}
\newtheorem{rem}[prop]{Remark}

\newtheorem{defi}[prop]{Definition}
\newtheorem{conj}[prop]{Conjecture/Question}

\def\begeq{\begin{equation}}
\def\endeq{\end{equation}}

\newcommand {\beg}{\begin{eqnarray*}}
\newcommand {\ee}{\end{eqnarray*}}

\def\and{\quad{\rm and}\quad}

\def\and{\quad{\rm and}\quad}

\let\lra=\longrightarrow

\def\mapright\#1{\,\smash{\mathop{\lra}\limits^{\#1}}\,}

\begin{document}
\bibliographystyle{plain}
\title{Space of K\"ahler metrics (IV)----On the lower bound of the K-energy }
\author{Xiuxiong Chen\footnote{The author is partially supported by NSF grant. }\\ Department of Mathematics\\
University of Wisconsin,  Madison}
\date{Sept. 12st, 2008}
\maketitle
\tableofcontents
\section{Introduction}
This is a continuation of \cite{chen05}.
In this paper, we want to partially confirm an old  conjecture by Donaldson \cite{Dona96}.
\begin{theo}  If there exists a  degenerating geodesic ray in the space of K\"ahler potentials tamed by a
bounded ambient geometry, then either  there is no constant scalar curvature (cscK) metric in this K\"ahler class, or there is a holomorphic line  consisting of cscK metrics only which is parallel to this degenerating geodesic ray.
\end{theo}

A complete resolution of Donaldson's conjecture is still out of reach at the present stage.  For detailed references in this subject, readers are encouraged to read \cite{chen991}\cite{chen992}\cite{chen05}.
A geodesic ray in the space of K\"ahler potential is called {\bf degenerated} ray if the K energy functional is non-increase along this ray. Definitions of the K-energy functional and degenerating geodesic rays  will be given in Section 2. We will also discuss more about various aspects of this theorem.  \\

Going back to Calabi's original view when he introduced  the extremal K\"ahler metric problem in \cite{calabi82}, one memorable feature is that the daring idealistic pursuit of K\"ahler geometry mix with
 hard  challenge in geometric  analysis.   While he painted a beautiful picture for us by promising that every K\"ahler class should admit a {\it best metric}---extremal K\"ahler metric in his sense, the daunting challenge  of solving it is close to  being``scaring'' since the Euler-Lagrange equation  is a fully nonlinear 6th order PDE.    In   \cite{calabi85}, E. Calabi showed us that the extremal K\"ahler metric always minimizes 
the Calabi energy locally.    When he did this work, we believe that he is only 60 years old and in the calculation of the second variation of the Calabi energy, he only use 16 derivatives to arrive at 
a positive sign with easy.  Since the day I became his student more than 20 years ago, I have been wondering all the time if there is the word ``fear'' in his mind.  20 years after his work,  we now know ( \cite{Dona052} in algebraic manifold  and \cite{chen05} in general K\"ahler manifold) that the extremal K\"ahler metric is also a global minimum of the Calabi energy in the underlying K\"ahler class.   Following \cite{chen991},   \cite{Dona05}  \cite{chentian005},  we know  that the cscK metric is a global minimizer
of the K-energy functional as well.   A natural question is:
\begin{conj} In any K\"ahler manifold $(M, [\omega])$, is $\displaystyle \inf_{\varphi \in [\omega]}\; Ca(\omega_{\varphi}) > 0 $ equivalent to $
 \displaystyle \inf_{\varphi \in [\omega]}\; \bE_{\omega}(\omega_{\varphi}) = -\infty$?
\end{conj}

One side of this question is already raised in \cite{chen05}: if $ \displaystyle \inf_{\varphi \in [\omega]}\; Ca(\omega_{\varphi}) = 0$, then the K-energy  bounds uniformly below in $[\omega]?\;$ According to \cite{chen05},   if $ \displaystyle \inf_{\varphi \in [\omega]}\; Ca(\omega_{\varphi}) = 0, $ then the $\yen$ invariant of all smooth degenerating geodesic rays vanishes. It follows that $(M, [\omega])$ is geodesic semi-stable with respect to all smooth geodesic rays.  Alternatively, in algebraic setting, by Donaldson \cite{Dona052}, the underlying polarization must
be Semi-K stable in the sense that generalized Donaldson-Futaki invariant 
must be non-positive.  In particular, the Calabi-Futaki invariant must vanish in $(M,[\omega]).\;$  One intriguing question is if  Semi-K-Stability or the vanishing of the Calabi-Futaki invariant alone is sufficient to imply the existence
of  a uniform lower bound on the K-energy functional in $(M,[\omega])$?   At least the answer to the second part of this question is negative, thanks to a beautiful example constructed by  S. K. Donaldson  \cite{Dona051}. Essentially by blowing-up $\CC P^{2}$ at various symmetric positions, Donaldson constructed a toric surface  where the Calabi-Futaki invariant
vanishes,   but the  K-energy functional approaches to $-\infty$ along at least one  degeneration.  In the same example,   he already showed us  that, along this degeneration, 
the Donaldson-Futaki invariant  is strictly negative.  Consequently, it follows by another theorem of Donaldson \cite{Dona052} in algebraic case (or \cite{chen05} in general), that the infimum of the Calabi energy in this K\"ahler manifold must be strictly positive.    Thus,  the question raised in \cite{chen05} (or at the beginning of this paragraph) survives this subtle example of Donaldson.  Nonetheless, this is a highly interesting example and  the author believes  that more can be learnt from
this example.  For instance, we may also check that  
the $\yen$ invariant of the corresponding geodesic ray (parallel to this degeneration) must be strictly
negative (c.f. \cite{chen05}), etc.  \\

This conjecture/question, as it states,  is rather general, abstract and elusive in many aspects.  Perhaps we can bring it down to earth by breaking it into several steps. 

\begin{conj} Is any of the following statements implies that the underlying K\"ahler class is  Semi-K-stable (when applicable) and the K-energy functional in this class is uniformly bounded from below?
\begin{enumerate}
\item there exists a cscK metric;
\item there exists a degenerating ray such that the K-energy  is bounded from below and the infimum of the Calabi energy along this ray is $0;\:$
\item  there exists a  {\bf critical} sequence of K\"ahler metrics;
\item the infimum of the Calabi energy over the space of K\"ahler potentials is $0.\;$

\end{enumerate}

\end{conj}
One can of course ask similar questions  by replacing geodesic ray with test configuration. 
Here we called a sequence of K\"ahler metrics $\{\varphi_{i}\} (i\in \NN)$  {\bf critical} when 
a) $\displaystyle \sup_{i\in \NN }\; |\bE_{\omega}|(\varphi_{i}) < \infty$; b) $ \displaystyle \inf_{i\in\NN}\;Ca(\omega_{\varphi_{i}}) = 0.\;$ 
Here the first statement is well-known by now (c.f. \cite{chen991},\cite{chentian005} and \cite{Dona05} ). 
The  4th statement is exactly one side of Conjecture 1.2.   For any convex function in finite dimensional
manifold with non-positive curvature,  the  second statement is correct, while it is not clear about the 3rd statement. There is a counter example to the 4th statement in finite dimensional manifold with non-positive curvature.  However, it is not clear if one can ``realize" such an example in the space of K\"ahler metrics.  Obviously, it is intriguing either to find such a counter example or prove the part 4 of Conjecture 1.3(readers are encouraged to read Section 3 for
more discussions on this direction).
Note that Conjecture/Question 1.3.2  is a natural extension of well known results  (c.f. \cite{chen991}, \cite{Dona05} and \cite{chentian005}).   We can prove this conjecture in some special case

\begin{theo} In toric variety, Conjecture 1.3.2 holds among smooth toric invariant K\"ahler metrics.
\end{theo}

   In fact, this theorem  holds whenever the geodesic conjecture of Donaldson \cite{Dona96} holds.  In general, we can answer Conjecture 1.3.2 positively once we strengthen the condition
slightly.
 
\begin{theo} If there exists a destabilized geodesic ray tamed by a bounded ambient geometry such that the K-energy functional is bounded from below along that particular ray;  and if the infimum of the Calabi energy is $0$ in this ray, then the K energy functional has a uniform  lower bound in the underlying K\"ahler class.   Moreover, this K\"ahler class is Semi-K-stable (whenever applicable).

\end{theo}

\begin{rem}  This is a surprising theorem since we only bound the K-energy  along one direction, while
the space of K\"ahler potential is infinite dimensional. The second condition on the ray is really crucial.  \\

\end{rem}

As an application, we prove the following theorem.   In some way, this  will give us  plenty of examples where the K-energy  is bounded from below in a  K\"ahler class which doesn't admit a  cscK metric. \\

\begin{theo}  For any destabilized simple test configuration such that the central fibre admits a cscK metric,   then the K-energy functional in the K\"ahler class defined by nearby fibre is bounded from below uniformly.  In particular, this K\"ahler class of nearby fibre is  Semi-K-stable(whenever applicable).
\end{theo}

 Note that  a K\"ahler class is Semi-K-stable if, for all test configurations,  the generalized Calabi-Futaki invariant at the central fibre has a preferred sign.  However, our theorem asserts that, as long as we have one simple test configuration where the central fibre is ``good " in some sense (admitting cscK metric), it is enough
 to conclude that the underlying polarization is Semi-K-stable!  This theorem might  be extended to more general setting. We delay more discussions on possible generalizations to the last section.   \\
 
     We may view this theorem from a different angle now.  The existence of a uniform lower bound on K-energy  had been always associated with the existence of cscK metric in the same K\"ahler class.  There is  no explicit example \footnote{no such example known to this author, but perhaps some experts knows.} where a lower bound on the K energy exists but no cscK metric in this K\"ahler class.  It is then worthwhile to review MuKai-Umermura 3 folds here.  Roughly speaking, it is a compactification of $PSL(2;\CC)/\Gamma$, where $\Gamma$ is
the icosahedra group(a subgroup of $SO(3)$ of order 60).   This gives a large family of 3-folds where rich geometry lives!  According to G. Tian \cite{tian97},  one of these three folds, called it $\M_1$ for easy of notations, admits  no KE metric in its canonical K\"ahler class.   This is a famous example where the classical Calabi conjecture on the existence of KE metric fails while the Futaki invariant vanishes.  Recently,
  Donaldson \cite{Dona081} showed  that, there is another 3-fold (denote as $\M_0$ ) in this
  family of 3-folds which lies in the closure of the orbit of $\M_1.\;$ Donaldson explicitly calculate  Tian's $\alpha$-invariant in $\M_0$ and he found it  to be  $ {5\over 6} > {3\over {3+1}}.\;$  By a theorem of G. Tian,  he then concludes that 
  there is a KE metric in this 3-fold $\M_0.\;$  \\
  
  Examining this family of complex 3-folds more closely, we can understand a lot more about
  deformation theory of the underlying complex structures.   
  First of all, we notice that $\M_0$ destabilizes $\M_1$ in the sense of Definition 1.7 of \cite{chen05}. One can even construct a simple test configuration  on MuKai-Umermura 3-fold, where the K\"ahler structure at nearby fibre is biholomorphic to $\M_1$ while the central fibre is biholomorphic to  $\M_0.\;$ By the preceding theorem, we know that the K energy functional in $\M_1$ must have a uniform lower bound in its canonical K\"ahler class although it admits no KE metric in its canonical class.  In fact, more can be said about geometry of MuKai-Umermura 3-folds.  For instance,  the author and his student \cite{chensong081} can prove that,  both K\"ahler Ricci flow and the Calabi flow initiated from some metric in canonical class in $\M_1$  will converge smoothly to the cscK metric at $\M_0.\;$ This certainly fits in well with the picture described by Theorem 1.7.   \\
  
  A result similar to Theorem 1.7 will hold if a K\"ahler manifold $(M, [\omega], J)$ is destabilized by another K\"ahler manifold $(M, [\omega]', J')$ where the destabilizer admits a cscK metric.   Definition of
  destabilizer is introduced in \cite{chen05}.  For Theorem 1.7 to be hold in this case, we need to modify the definition 1.7 \cite{chen05} of destabilization slightly (according to Theorem 1.5 above). In a sequeal paper of this one, the author and his student will  discuss alternative
  proof of results in algebraic setting \cite{chensong082}. \\
  
  \noindent {\bf Organization} In Section  2, we give a brief account of K\"ahler geometry including a
  proof of subharmonicity of the K energy over smooth  solution on geodesic equation of disc version. 
  In Section 3, we give a proof of Theorem 1.4 and 1.7. In some sense, Theorem 1.7 is a corollary or application of Theorem 1.5.  However, the proof of Theorem 1.7 at this section is independent of Theorem 1.5.  In Section 4, we give a proof to Theorem 1.5.  In Section 5, we prove   Theorem 1.1.\\
 
    \noindent {\bf Acknowledgment} While the author obtained main results in this paper much earlier, the main writing is done during his visit to Tokyo Institute of Technology and University of Science and
Technology of China this summer.   The author is very grateful to the wonderful hospitality provided by Tokyo-Tech and
USTC respectively. The author would like to thank Professor A. Futaki for bringing Donaldson's example to  his attentions. He also wants to thanks colleague(s), graduate student(s) at Tokyo-Tech for the kindness they have shown during his stay at Tokyo-Tech.  The author wishes to thank his student S. Sun for help with examples.\\

The author lists a number of problems in this paper. He wishes to acknowledge that some of these problems may have been known to various experts in the field.  It is hope that this (by putting them together) is helpful to younger generation K\"ahler geometers. 
\section{Brief outline of geometry in the Space of K\"ahler potentials.}
\subsection{Quick introduction of K\"ahler geometry}
Let $\omega$ be a fixed K\"ahler metric
on $M$. In a local holomorphic coordinate, $\omega$ can be expressed as
\[
\omega = {\sqrt{-1}\over 2}\;  g_{\alpha\b \beta} \;d w^\alpha \wedge
d\,w^{\b \beta}  >  0.
\]

The scalar curvature can be defined as
\[
R  = - g^{\alpha\b \beta }\;  {{\p^2 \log \det \left(g_{i\b
j}\right) }\over {\p w^\alpha \p w^{\b \beta}}}.
\]
The so-called Calabi energy is
\begin{equation}
Ca(\omega) = \int_M \; (R(\omega) - \underline R)^2 \omega^n,\;\qquad {\rm where}\qquad \underline R = {{[C_{1}(M)] \cdot [\omega]^{[n-1]}}\over [\omega]^{[n]}}.
\label{eq:calabienergy}
\end{equation}
 According to Calabi 
\cite{calabi82}\;\cite{calabi85}, a K\"ahler metric is called extremal if the
complex gradient vector field
\begin{equation}
\cX_c = g^{\alpha\b \beta} {{\p R}\over {\p w^{\b \beta}}}
{\p\over {\p w^\alpha}} \label{eq:extremalvectorfield}
\end{equation}
is a holomorphic vector field.  \\

If $X$ is a holomorphic vector field, then for any K\"ahler potential
$\varphi$ we can define $\theta_X$ up to some additive constants by
\begin{equation}
L_X \omega_\varphi = \sqrt{-1} \p \b \p\; \theta_X (\varphi).
\label{eq:liederivatives}
\end{equation}
Then, the well known Calabi-Futaki invariant
\cite{Futaki83}\;\cite{calabi85} is 
\begin{equation}
\cF_X([\omega]) =  \displaystyle
\int_M\; \theta_X (\varphi) \cdot (\underline{R} - R(\varphi))
\;\omega_\varphi^n. \label{eq:calabifutakiinvariant}
\end{equation} Note that this is a Lie algebra
character which depends on the K\"ahler class only.

\subsection{Weil-Petersson type metric by Mabuchi}
Denote the space of smooth K\"ahler potentials as
\[
{\cH}= \{ \varphi \mid \omega_{\varphi} = \omega + \b \p
\p \varphi > 0,\;{\rm on} \; M\}/ \sim,
\]
where $\varphi_1 \sim \varphi_2$ if and only if
$\varphi_1=\varphi_2+ c$ for some constant $c$. A tangent
vector in $T_\varphi\cH$ is just a function $\psi$ such that
$$\int_M \psi\omega_\varphi^n =0.$$ Its norm in the $L^2$-metric
on $\cH$ is given by (cf. \cite{Ma87})
\[
\|\psi\|^2_{\varphi} =\int_{M}\psi^2\;\omega^n_\varphi.
\]
This metric was subsequently re-defined in \cite{Semmes92} and
\cite{Dona96}.    In all three papers, \cite{Ma87}\cite{Semmes92} and \cite{Dona96},  the authors
defined  this Weil-Petersson type metric from various points of view and proved formally that this infinite
dimensional space has non-positive curvature.   Using this definition,  we can define a distance
function in $\cH$:  For any two K\"ahler potentials $\varphi_0,
\varphi_1\in \cH$, let $d(\varphi_0,\varphi_1)$ be the infimum of
the length of all possible curves in $\cH$ that connect $\varphi_0$
with $\varphi_1.\;$  \\

A straightforward computation shows that a geodesic path
$\varphi:[0,1]\rightarrow \cH$ of this $L^2$ metric must satisfy
the following equation
\[
\varphi''(t) - {g_\varphi}^{\alpha \b \beta} {{\p^2 \varphi}\over
{\p t \p w^\alpha}}  {{\p^2 \varphi}\over {\p t \p w^{\b \beta}}}
= 0.
\]
where
\[
g_{\varphi,\alpha\b \beta} =  g_{\alpha \b \beta} + {{\p^2
\varphi}\over {\p w^\alpha\p w^{\b \beta}}} > 0.
\]
 
According to S. Semmes \cite{Semmes92},   a smooth path $\{\phi(t), t\in [0,1] \} \subset \cH$ satisfies the
geodesic equation if and only if the function $\phi$ on $[0,1]\times
S^1\times M$ satisfies the homogeneous complex Monge-Ampere equation
\begin{equation}
\label{eq:hcma0} (\pi_2^* \omega + \partial \overline{\partial
}\phi)^{n+1} \;\; =\;\; 0, \qquad {\rm on} \; \Sigma \times M,
\end{equation}
where $\Sigma = [0,1] \times S^1$ and $\pi_2: \Sigma\times M\mapsto
M$ is the projection.  In fact, one can consider (\ref{eq:hcma0})
over any Riemann surface $\Sigma$ with boundary condition
$\phi =\phi_0$ along $\partial \Sigma$, where $\phi_0$ is a smooth
function on $\partial \Sigma \times M$ such that $\phi_0(z,\cdot)\in
\cH$ for each $z\in \partial \Sigma$.\footnote{We often regard
$\phi_0$ as a smooth map from $\p \Sigma$ into $\cH$.}  The equation (\ref{eq:hcma0}) can be regarded as
the infinite dimensional version of the WZW equation for maps from
$\Sigma$ into $\cH$ (cf. \cite{Dona96}).\footnote{The original WZW
equation is for maps from a Riemann surface into a Lie group.}
\\

Next we introduce a  well known  functional in $\cH$
here.    The 
 K-energy functional
(introduced by T. Mabuchi) is defined as a closed form $d\,\bE$.
Namely, for any $\psi \in T_\varphi \cH$, we have
\begin{equation}
(d\,\bE, \psi)_\varphi = \int_M\; \psi\cdot (\underline{R} -
R(\varphi)) \;\omega_\varphi^n.\label{eq:kenergyform}
\end{equation}
Note that for any holomorphic vector field, we have
\[
\cF_X([\omega]) = (d\,\bE, \theta_X)_\varphi.
\]

\subsection{Geodesic rays and geodesic stability} 
  In this subsection, we collect some definitions, theorems from \cite{chen05} which are useful
  in this paper.

\begin{defi}\label{def:specialgeodesicray} \cite{chen05} A smooth geodesic ray
$\rho(t) (t \in [0,\infty))$ is called special if it is one of the following types:

\begin{enumerate}
\item {\bf effective} if  the Calabi energy of $\omega_{\rho}$ in $M$ is dominated by ${\epsilon\cdot  t^2}$ for any $\epsilon > 0$ as $t \rightarrow \infty.\;$
\item{\bf  normal}  if the curvature of $\omega_{\rho} $ in $M$ is uniformly bounded  for $t\in [0,\infty).\;$
\item {\bf bounded geometry} if  $(M, \omega_{\rho(t)}) (t\in [0,\infty)) $  has uniform bounds on  curvature and uniform positive lower bounds on the injective radius. 
\end{enumerate}
\end{defi}

\begin{defi} {\bf Bounded ambient geometry} A K\"ahler metric $h = \pi_2^* \omega_0 + i \p \b \p \b \rho$ 
in  $([0,\infty)\times S^1) \times M$ is said to have bounded ambient geometry if
\begin{enumerate}
\item  it has a uniform bound on its curvature;
\item  $([0,T]\times S^1 \times M, h)$ has a uniform lower bound on injectivity radius
 and the bound is independent of $T\rightarrow \infty;\;$
\item The vector length $|{\p\over {\p t}}|_h$ has a uniform upper bound. 
 \end{enumerate}

Let $h$ be a bounded ambient  K\"ahler  metric in $[0,\infty)\times S^{1}\times  M.\;$
Suppose that its corresponding K\"ahler form $\tilde{\omega}$ is given
by\footnote{Here $\pi_2: ([0,\infty)\times S^1 ) \times M \rightarrow M$ be the natural projection map.}
\begin{equation}
  \pi_2^* \omega_0 + \displaystyle \sum_{i,j=1}^n\; {{\p^2 \b \rho}\over
  {\p w^i \p w^{\b j}}}  d\, w^i d\, w^{\b j}  + 2 Re\left(\displaystyle \sum_{i=1}^n\; {{\p^2 \b \rho}\over
  {\p w^i \p \b z}} d\, w^i d\, \bar z \right) +  {{\p^2  \b \rho}\over {\p z \p \b z}}
  d\, z\; d\, \b z.
\label{eq:wrappedproductmetric0} \end{equation} Here $z = t +
\sqrt{-1} \theta.\;$ In other words
\begin{equation}
\tilde{\omega} = \pi_2^*\omega_0 + \sqrt{-1} \p \b \p \b \rho .
\label{eq:wrappedproductmetric1}
\end{equation}
In the remainder of this section, we will use $h$ to denote any K\"ahler  metric in $[0,\infty)\times S^{1}\times M$  with bounded ambient geometry.

\end{defi} 
\begin{defi} {\bf Tamed by a bounded ambient geometry}     A smooth geodesic ray $ \rho:[0,\infty) \rightarrow \cH$  is said to be tamed by a bounded ambient geometry $h$, if there is a uniform bound for the relative potential $\rho - \bar{\rho}$.
 \end{defi}

For most purposes, ``weakly tamed by a bounded ambient geometry weakly'' is sufficient.
In fact, any normal geodesic ray is expected to be tamed by some bounded ambient
geometry, at least when it has bounded geometry.    Any effective geodesic ray is expected to be tamed
by some bounded ambient geometry.  Since the space of K\"ahler potential is non-positively curved in
the sense of Alexandrov \cite{chen992}, it is natural to introduce the notion of parallelism between two geodesic rays.

\begin{defi} \cite{chen05} \label{def:geodesicparallel} Any two smooth geodesic rays $\rho_{1},  \rho_2:
[0,\infty)\rightarrow \cH,\;$  are called parallel if there
exists a constant $ C$ such that
\[
  \displaystyle \sup_{t \in [0,\infty)}\; |\rho_1(t)-\rho_2(t)|  \leq C.
\]
\end{defi}

For any smooth geodesic ray, there is some geodesic ray (perhaps in the weak form) which
parallels to it as in Definition 2.5, but not necessary relatively  $C^{1,1}$ (with respect to the original rays).  In \cite{chen05}, we made some attempt in this direction and we quote our modest result here.

\begin{theo} \cite{chen05}\label{th:geodesicDegenerate} If there exists a  smooth geodesic ray  $\rho: [0,\infty) \rightarrow \cH$ which is tamed by a bounded ambient geometry,
then for any K\"ahler potential $\varphi_0\in \cH$,  there exists a relative $C^{1,1}$ geodesic ray $\varphi(t)$ initiated from
$\varphi_0$ and parallel to $\rho.\;$ 
 \end{theo}

Since the K energy is well defined for $C^{1,1}$ K\"ahler potential, so we defined an invariant
as 

\begin{defi} \label{defi:geodesicinvariant} For every relative $C^{1,1}$ geodesic ray $\varphi:  [0,\infty) \rightarrow \cH$,  we
can define an invariant as
\begin{equation} \yen(\varphi) =  \displaystyle \liminf_{t \rightarrow
\infty}  (\bE(\varphi(t) - \bE(\varphi(t-1))).\label{eq:invariant2}
\end{equation}
\end{defi}

Clearly, in the case of smooth geodesic ray, this limit agrees with the invariant introduced in \cite{chen05}: 
\begin{defi} \label{defi:geodesicinvariant} \cite{chen05}For every smooth geodesic ray $\rho:  [0,\infty)\rightarrow \cH$,  we
can define an invariant as
\begin{equation} \yen(\rho) =  \displaystyle \lim_{t \rightarrow
\infty} \displaystyle \int_M\; {{\p \rho(t)}\over {\p t}}
(\underline{R} - R(\rho)) \omega_{\rho}^n.
\label{eq:invariant1}
\end{equation}
\end{defi}

Now we introduce notions of geodesic stability which already appear in \cite{chen05}.

\begin{defi}\label{def:geodesicdestable} \cite{chen05}
A smooth geodesic ray $\rho: [0,\infty) \rightarrow \cH$ is called
stable (resp; semi-stable) if $\yen(\rho)>0$ (resp: $\geq 0$).   It is called a destabilizer for
$\cH$ if $\yen(\rho) < 0.\;$ \end{defi}

Following the approach used in the algebraic case,  we define (cf. \cite{Dona96}):
\begin{defi}\label{def:geodesicstable} \cite{chen05} A K\"ahler  manifold is called geodesically stable if
there is no  destabilizing smooth geodesic ray.  It is called
weakly geodesically stable if the invariant $\yen$ is always
non-negative for every smooth  geodesic ray.
\end{defi}

\subsection{Subhamonicity of the K energy over WZW solutions}

Suppose that $\phi$ is a
$C^{1,1}$ solution of (\ref{eq:hcma0}), we denote by $\cR_\phi$
the set of all $(z,x)\in \Sigma\times M$ near which $\phi$ is
smooth and $\omega_{\phi(z',\cdot)}=\omega
+\sqrt{-1}\partial\overline{\partial} \phi(z',\cdot)$ is a
K\"ahler metric. We may regard $\cR_\phi$ as the regular set of
$\phi$. It is open, but {\it a priori}, it may be empty. We have a
distribution $\cD_\phi\subset T(\Sigma\times M)$ over $\cR_\phi$:
\begin{eqnarray}
\label{eq:cDovercR} \cD_\phi|_{(z,x)} = \{ v \in T_z\Sigma\times
T_xM~|~i_v \left (\pi^*_2\omega + \sqrt{-1} \p\overline{\p}
\phi\right ) =0\},\qquad \forall\;(z,x)\in \cR_\phi.
\end{eqnarray}
Here $i_v$ denotes the interior product. Since the form is closed,
$\cD_\phi$ is integrable when $\phi$ is smooth. We say that $\cR_\phi$ is saturated in
$\cV\subset \Sigma\times M$ if every maximal integral sub-manifold
of $\cD_\phi$ in $\cR_\phi\cap \cV$ is a closed disk  in $\cV$.
On any product manifold, we may write any vector in
$\cD_\phi$ as
\begin{equation}
{\p \over {\p z}} + X \in \cD_\phi\mid_{(z,x)},\qquad {\rm
where}\;\; X \in T^{1,0}_x M. \label{def:leafvector}
\end{equation}
\begin{defi}
\label{def:partiallysmooth} A solution $\phi$ of (\ref{eq:hcma0})
is called partially smooth if it is $C^{1,1}$-bounded on
$\Sigma\times M$ and $\cR_\phi$ is open and saturated in
$\Sigma\times M$, but dense in $\p \Sigma \times M$, such that the
varying volume form $\omega_{\phi(z,\cdot)}^n$ extends to a
continuous $(n,n)$ form on $\Sigma^0\times M$, where
$\Sigma^0=(\Sigma\backslash \partial \Sigma)$.
\end{defi}

Clearly, if $\phi$ is a partially smooth solution, then its
regular set $\cR_\phi$ consists of all points where the vertical
volume form $\omega_{\phi(z,\cdot)}^n$ is positive in
$\Sigma\times M$.



\begin{defi}
\label{def:almostsmooth} We say that a solution $\phi$ of
(\ref{eq:hcma0}) is almost smooth if
\begin{enumerate} \item it is partially smooth,
\item The distribution ${\cal D}_\phi$ extends to a continuous
distribution in a saturated set $\tilde \cV\subset \Sigma \times
M$, such that the complement $\tilde \cS_\phi$ of $\tilde \cV$ is
codimension at least 2
 and $\phi$ is
$C^1$ continuous on $\tilde \cV$. The set $\tilde \cS_\phi$ is
referred to as the singular set of $\phi$. \item The leaf vector
field $X$ is uniformly bounded in $\cD_\phi.\;$
\end{enumerate}
\end{defi}

A smooth solution is certainly an almost smooth solution of
(\ref{eq:hcma0}). If the boundary values of a sequence of almost smooth solutions converge in some $C^{k,\beta}$ topology, then the sequence
converges to a partially smooth solution in the $C^{1,\beta}$-topology
($ k > 2, \;0< \beta < 1$).

\begin{theo}\cite{chentian005}
\label{th:almostsmooth} Suppose that $\Sigma$ is a unit disc. For
any $C^{k,\alpha}$ map $\phi_0: \p \Sigma \rightarrow \cH$ ($k \ge
2$, $0< \alpha < 1$) and for any $\epsilon
> 0$, there exists a $\phi_\epsilon: \p\Sigma\rightarrow \cH$
in the $\epsilon$-neighborhood of $\phi_0$ in
$C^{k,\alpha}(\Sigma\times M)$-norm, such that (\ref{eq:hcma0})
has an almost smooth solution with boundary value $\phi_\epsilon$.
\end{theo}

An almost smooth solution of eq. \ref{eq:hcma0} has uniform
$C^{1,1}$ bounds and is smooth almost everywhere.   The K-energy functional is well defined
for this family of K\"ahler potentials. 

\begin{theo} \cite{chentian005}\label{appl:kenergyissubharmonic} Suppose that $\phi: \Sigma \rightarrow \overline{{\cal H}_\omega}$
is an almost smooth solution described as in Definition
\ref{def:almostsmooth}. Then the induced K-energy function $\bE:
\Sigma\rightarrow \RR$ (by $\bE(z) = \bE(\phi(z,\cdot))$) is
weakly sub-harmonic and $C^1$ continuous (up to the boundary).
More precisely,
\[
{{\p^2}\over {\p z\p\b z}}  \bE(\phi(z,\cdot)) = \int_{\pi\circ
ev(z,\cU_{\phi_0})} | {\cal D} {{\p \phi}\over {\p \b
z}}|_{{\omega_{\phi}}}^2\, {\omega_{\phi}}^n\;\geq 0, \qquad
\forall\; z \in \Sigma^0
\] holds in $\Sigma^0$ in the weak sense. On $\p \Sigma$, we have
\[
 \displaystyle \int_{\p
\Sigma}{{\p \bE}\over{\p\,{\bf n}}} (\phi) d s
= \displaystyle \int_{\pi\circ
ev(z,\cU_{\phi_0})} | {\cal D} {{\p \phi}\over {\p \b
z}}|_{{\omega_{\phi}}}^2\, {\omega_{\phi}}^n\, d s,
\] where  $d s$
is the length element of $\p \Sigma,\;$ and $\bf n$ is the outward
pointing unit normal direction at $\p \Sigma.\;$

\end{theo}

To help readers to understand this theorem better, we will present a proof of this theorem in the case that the disc version geodesic solution $\phi$ is smooth.  Hence this is just a formal proof and the readers are encourage to read \cite{chentian005} for a vigorous proof of this fact.  Note that for any smooth path $\phi(t)$, we have
\[
{{d^2 \bE}\over {d\, t^2}} (\phi(t)) = - \displaystyle  \int_M\; ({{\p^2 \phi}\over {\p t^2}} - |\nabla {{\p \phi}\over {\p t}}|_\phi) (R_\phi- \b R) \omega_{\phi(t)}^n  +\displaystyle \int_M\; |{\cal D}  {{\p \phi}\over {\p t}}|^2_\phi \omega_\phi^n.
\]
Note that for a disc version geodesic, we have
\[
z = t+ \sqrt{-1} s, \qquad {\rm and}\; \triangle_\Sigma =  {{\p^2 }\over {\p t^2}} +    {{\p^2 }\over {\p s^2}}.
\]
Thus,
\[
\begin{array}{lcl} \Delta_\Sigma \bE & = &  ({{\p^2 }\over {\p t^2}} +    {{\p^2 }\over {\p s^2}}) \bE\\ & = &  - \displaystyle  \int_M\; ({{\p^2 \phi}\over {\p t^2}} - |\nabla {{\p \phi}\over {\p t}}|_\phi^2) (R_\phi- \b R) \omega_{\phi(t,s)}^n    - \displaystyle  \int_M\; ({{\p^2 \phi}\over {\p s^2}} - |\nabla {{\p \phi}\over {\p s}}|_\phi^2) (R_\phi- \b R) \omega_{\phi(t,s)}^n \\
&&\qquad +\displaystyle \int_M\; |{\cal D}  {{\p \phi}\over {\p t}}|^2_\phi \omega_\phi^n  +\displaystyle \int_M\; |{\cal D}  {{\p \phi}\over {\p s}}|^2_\phi \omega_\phi^n.
\end{array}
\]
Using the equation for disc version geodesics, we have
\[
\begin{array}{lcl} \Delta_\Sigma \phi  & = & g_\phi^{\alpha \b \beta} ({{\p \phi}\over {\p \b z}})_{\alpha} ({{\p \phi}\over {\p z}})_{\b \beta}  = g_\phi^{\alpha \b \beta} ({{\p \phi}\over {\p t}} - \sqrt{-1} {{\p \phi}\over {\p s}})_{\alpha} ({{\p \phi}\over {\p t}} +  \sqrt{-1} {{\p \phi}\over {\p s}})_{\b \beta} \\
& = & |\nabla {{\p \phi}\over {\p t}}|^2_\phi + |\nabla {{\p \phi}\over {\p s}}|^2_\phi + \sqrt{-1} g_\phi^{\alpha\b \beta}  \left( ({{\p \phi}\over {\p t}})_\alpha ({{\p \phi}\over {\p s}})_{\b \beta} -  ({{\p \phi}\over {\p s}})_\alpha ({{\p \phi}\over {\p t}})_{\b \beta} \right)\\
& = &  |\nabla {{\p \phi}\over {\p t}}|^2_\phi + |\nabla {{\p \phi}\over {\p s}}|^2_\phi  + \{{{\p \phi}\over {\p t}}, {{\p \phi}\over {\p s}} \}_\phi.
\end{array}
\]
The last term gives the Poisson bracket with respect to the symplectic form $\omega_\phi$. Thus,
we have
\[
\begin{array}{lcl} \Delta_\Sigma \bE & = & - \displaystyle  \int_M\; ({{\p^2 \phi}\over {\p t^2}}  + {{\p^2 \phi}\over {\p s^2}}- |\nabla {{\p \phi}\over {\p t}}|_\phi^2 -  |\nabla {{\p \phi}\over {\p t}}|_\phi^2) (R_\phi- \b R) \omega_{\phi(t,s)}^n    \\
&&\qquad +\displaystyle \int_M\; |{\cal D}  {{\p \phi}\over {\p t}}|^2_\phi \omega_\phi^n  +\displaystyle \int_M\; |{\cal D}  {{\p \phi}\over {\p s}}|^2_\phi \omega_\phi^n\\
& = & -( \{{{\p \phi}\over {\p t}}, {{\p \phi}\over {\p s}} \}_\phi, R_\phi - \bar R)_\phi +\displaystyle \int_M\; |{\cal D}  {{\p \phi}\over {\p s}}|^2_\phi \omega_\phi^n  + \displaystyle \int_M\; |{\cal D}  {{\p \phi}\over {\p t}}|^2_\phi \omega_\phi^n.
\end{array}
\]
The first term in the last line may give us some trouble.  However, one notices that
\[
\begin{array}{lcl}  \displaystyle \int_M\; |{\cal D} {{\p \phi}\over {\p \b z}} |^2_\phi \omega_\phi^n & = & ({\cal D} ( {{\p \phi}\over {\p t}} - \sqrt{-1} {{\p \phi}\over {\p s}}), {\cal D} ({{\p \phi}\over {\p t}} - \sqrt{-1} {{\p \phi}\over {\p s}}))_\phi\\
& = &  \displaystyle \int_M\; |{\cal D}  {{\p \phi}\over {\p t}}|^2_\phi \omega_\phi^n + \displaystyle \int_M\; |{\cal D}  {{\p \phi}\over {\p s}}|^2_\phi \omega_\phi^n + \sqrt{-1}  ({\cal D}  {{\p \phi}\over {\p t}}, {\cal D}  {{\p \phi}\over {\p s}})_\phi - \sqrt{-1}  ({\cal D} {{\p \phi}\over {\p s}}, {\cal D} {{\p \phi}\over {\p t}})_\phi\\
& = &  \displaystyle \int_M\; |{\cal D}  {{\p \phi}\over {\p t}}|^2_\phi \omega_\phi^n + \displaystyle \int_M\; |{\cal D}  {{\p \phi}\over {\p s}}|^2_\phi \omega_\phi^n + \sqrt{-1}  ( (\b {\cal D} {\cal D} -{\cal D}\b D)  {{\p \phi}\over {\p t}},  {{\p \phi}\over {\p s}})_\phi\\
& = &  \displaystyle \int_M\; |{\cal D}  {{\p \phi}\over {\p t}}|^2_\phi \omega_\phi^n + \displaystyle \int_M\; |{\cal D}  {{\p \phi}\over {\p s}}|^2_\phi \omega_\phi^n + \sqrt{-1} \int_M\; g_\phi^{\alpha \b \beta} ( R_\alpha ( {{\p \phi}\over {\p t}})_{\b \beta} -  R_{\b \beta} ( {{\p \phi}\over {\p t}})_{\alpha}) {{\p \phi}\over {\p s}} \omega_\phi^n\\
& = &\displaystyle \int_M\; |{\cal D}  {{\p \phi}\over {\p t}}|^2_\phi \omega_\phi^n + \displaystyle \int_M\; |{\cal D}  {{\p \phi}\over {\p s}}|^2_\phi \omega_\phi^n +  (\{ R,  {{\p \phi}\over {\p t}}\}_\phi, {{\p \phi}\over {\p s}})_\phi\\
& = & \displaystyle \int_M\; |{\cal D}  {{\p \phi}\over {\p t}}|^2_\phi \omega_\phi^n + \displaystyle \int_M\; |{\cal D}  {{\p \phi}\over {\p s}}|^2_\phi \omega_\phi^n +  (\{{{\p \phi}\over {\p t}}, {{\p \phi}\over {\p s}}\}_\phi, R -\bar R)_\phi.
\end{array}
\]
Now, plugging this into $\triangle_\Sigma \bE$, we have
\[
 \Delta_\Sigma \bE_\omega =  \displaystyle \int_M\; |{\cal D} {{\p \phi}\over {\p \b z}} |^2_\phi \omega_\phi^n  \geq 0.
\]
Thus, we prove subharmonicity for any smooth disc version geodesic. \\

\section{Proof of Theorem 1.4 and 1.7}
Let us state a lemma in elementary calculus first.
\begin{lem} Suppose $f$ is a convex function in $\RR^{2}$ such that $f$ is non-increasing and uniformly bounded from below alone
one half line. If in additionally we assume $\liminf |\nabla f|$ alone this half line is  $0$, then $f$ is bounded from below globally. 
\end{lem}
We remark that  this lemma hold in any manifold with non-positive curvature.  Following \cite{chen992}, the space of K\"ahler potentials is non-positive in the sense of Alexanderov.  In the case of toric variety, there is a unique smooth geodesic segment between any two invariant K\"ahler potentials.  The K energy is convex on smooth geodesic segment while the Calabi energy is precisely the norm of the
``gradient'' of the K energy functional in this infinite dimensional space. Thus,  the proof of Theorem 1.4 is exactly the same as the proof of this lemma.  Now we will give a proof to this lemma
only.  
\begin{proof} Suppose $P_{i} =(x_{i}, y_{i})$ is a sequence of points in this half line  in $\RR^{2}$ which diverges to $\infty\;$ such that \[
|\nabla f(x_{i}, y_{i}) | \rightarrow 0, \qquad {\rm and}\;\; f(x_{i}, y_{i}) \geq -C
\]
for some uniform constant $C> 0.\;$  Let $O=(0,0)$ be the original point of this half line.  Let $Q=(x,y)$ be any other point in $\RR^{2}.\;$  We want to argue that $f(Q)$ is bounded from below by a constant independent of $q.\;$\\

Consider the sequence of triangle $OP_{i}Q (i \in \NN).\;$ Set
\[
|OQ|=d,\qquad |OP_{i}| = l_{i} \rightarrow \infty, \qquad {\rm and}\;\; |P_{i}Q| = \tilde l_{i}
\]
Then,
\[
 l_{i} - d < \tilde l_{i} < l_{i} + d, \qquad \forall i \in \NN.
\]
Further more, when $i$ large enough, we have
\[
0 < \sin \angle OP_{i}Q  \leq {{10 d}\over l_{i}}.
\]
Let $u_{i}, v_{i}$ be the two unit  tangent vectors at $T_{P_{i}} \RR^{2}$ pointing to $Q, O$ respectively.
Then, \[
(u_{i}- v_{i}, u_{i} - v_{i}) = {{100 d^{2}}\over l_{i}^{2}}
\]
Since $P_{i}$ is in the same half line where the functional $f$ is decreasing, then
\[
 (\nabla f, v_{i})\mid_{P_{i}} \geq 0.
\]
Restricting the convex function $f$ in the line segment $QP_{i}$, we have
\[
\begin{array}{lcl} {{f(Q) - f(P_{i})} \over |P_{i}Q|} & \geq & ( \nabla f, u_{i})\mid_{P_{i}} \\
& = & (\nabla f, v_{i})\mid_{P_{i}}  + (\nabla f, u_{i} -  v_{i})\mid_{P_{i}} \\
& \geq &  -|\nabla f|_{ P_{i}}\cdot |u_{i} - v_{i}|_{P_{i}} = - |\nabla f|_{P_{i}} \cdot {{10 d}\over l_{i}}.
\end{array}
\]
Thus, we have
\[\begin{array}{lcl}
f(Q) & \geq & f(P_{i})  - |\nabla f|_{P_{i}}\cdot {{10 d}\over l_{i}} \cdot \tilde{ l}_{i} \\
& \geq & -C -  |\nabla f|_{P_{i}} \cdot 20 d.
\end{array}
\]
Now let $l_{i}\rightarrow \infty$,  we have
\[
f(Q) \geq -C.
\]
The lemma is then proved.

\end{proof}
From the proof, the second condition is absolutely necessary. In fact, we have the following counter example:
\[
  f (x, y) = e^{x} + y, \qquad (x,y) \in \RR^{2}.
\]
With respect to conjecture 1.3.4, it is false in finite dimensional as the following example suggested.
\[
f(x) = \left\{ \begin{array}{ll} {x^2 \over 2} - x, & \qquad {\rm if }\;\; x \leq 0,\\
- \ln (x+1) &\qquad {\rm if }\;\; x \geq 0. \end{array} \right.
\]
This is a convex function where
\[
\displaystyle \inf_{x \in \RR} \; |\nabla f | (x)  = 0 \qquad {\rm and}\qquad \displaystyle \inf_{x \in \RR} \;  f(x)  = -\infty.\]
On the other hand, the answer to conjecture 1.3.4 might still be correct because of the geometric structure attached to this problem.  \\

Let us first recall a theorem from \cite{chen05}.
\begin{theo} \cite{chen05} Let $\varphi_0, \varphi_1$ be two arbitrary smooth K\"ahler potentials in the same K\"ahler class.  Then, the following inequality holds
\begin{equation}
\bE (\varphi_1) - d(\varphi_0, \varphi_1) \cdot \sqrt{Ca(\varphi_1)} 
\leq \bE(\varphi_0).
\end{equation}
Here $d(\varphi_0, \varphi_1)$ is the geodesic distance in the space of K\"ahler potentials.
\end{theo}

 To prove Theorem 1.7, we only need a weak version of Theorem 1.5 whose proof is much simpler. \\

\noindent {\bf Theorem 1.5a. } {\it If there exists a destabilized geodesic ray such that the K-energy functional is bounded from below along this given ray and if the Calabi energy converges to $0$ in the speed of $o({1\over t^{2}})$ along the ray, then the K-energy functional has a uniform lower bound
 in the underlying
K\"ahler class.  (Here $t$ represents the distance along geodesic ray from the initial potential in the ray.) Moreover, this K\"ahler class is Semi-K-stable (when applicable).}\\

We give a proof to this weak version of Theorem 1.5a.

\begin{proof} Suppose $\rho: [0,\infty) \rightarrow \cH $ is a destabilized geodesic ray such that the K energy is uniformly bounded from below:

\[
   \bE (\rho) \geq -C, \qquad \forall t \in [0,\infty)
\]
where $C$ is some constant which change from line to line.   Our assumption means that

\begin{equation}
\lim_{t\rightarrow \infty}\; t \cdot Ca(\rho(t)) \leq C.
\end{equation}

 Let $\varphi_{0}$ be any K\"ahler potential in $\cal H$. For any $t\geq 1$, apply Theorem 3.2 to this case. By inequality (3.1),  we have
\[
\begin{array}{lcl}
E(\varphi_{0}) & \geq & E(\rho(t)) - d(\varphi(0),\rho(t)) \cdot \sqrt{Ca(\rho(t))}\\
& \geq & E(\rho(t)) - (d(\varphi(0), \rho(0)) + t) \cdot \sqrt{Ca(\rho(t))} \geq  -C.
\end{array}
\] 
It follows that the K-energy has a uniform lower bound in $[\omega].\;$
\end{proof}

 Before we give a proof of Theorem 1.7, we introduce an important proposition.
 
\begin{prop} \cite{chensong081} Given a simple test configuration
$\L\rightarrow(\M, J, \Omega)\rightarrow D$. Suppose there exists a smooth $S^1$ invariant solution of the homogeneous
 complex Monge-Amp\'{e}r equation $(\Omega+\sqrt{-1}\partial\overline{\partial}\Phi)^{n+1}=0$, and denote the induced
 geodesic ray by $\phi_t(t\in [0,\infty))$. If the central fiber has cscK metric, then
 \begin{enumerate}
 \item the first derivatives of the K energy decay exponentially. Consequently,  the
K-energy has a lower bound along the geodesic ray. 
\item there exists a geodesic ray from nearby fibre (tamed by this simple test configuration) where the Calabi energy decays exponentially as well.
\end{enumerate}
\end{prop}

We will delay the proof to the end of this section. We give a proof of Theorem 1.7 here.

\begin{proof} Suppose $\pi: \chi \rightarrow D$ with smooth central fiber.  According to Arezzo-Tian \cite{Are-Tian03}, we can prove the existence of geodesic ray parallel to this simple test configuration,  initiating from a nearby fibre (close enough) to  the central fibre. In \cite{song08},  S. Sun  gives an alternative proof to this beautiful existence problem.    According to Chen-Tang \cite{chentang07},  the $\yen$ invariant of the induced geodesic ray  is precisely the evaluation of the Calabi-Futaki invariant
in the central fibre.  If the Calabi Futaki invariant in the central fibre vanishes, then the $\yen$ invariant of induced geodesic ray must be $0$. In particular, this geodesic ray must be a destabilizing geodesic ray. \\

From the method of obtaining this geodesic ray,  one  can specify which limit metric of a geodesic ray constructed  from \cite{song08} will converges to (of course, up to a holomorphic transformation
in the central fibre). This is very crucial here: We need the geodesic ray converges to the cscK metric in the central fibre in the limit.  Then, we can have a ray whose K energy values bounded uniformly from below while the Calabi energy converges to $0.\;$ In other words, Prop. 3.3 may hold for this geodesic ray. Then, following Theorem 1.5a, we know that the K energy in the entire K\"ahler class must be bounded from below.
  \end{proof}

Now we are ready to prove Proposition 3.3.  
\begin{proof}
Let $\M$ be the total space of this simple test configuration $\M \rightarrow \triangle.\;$ Here $\triangle$
is a unit disc in $\CC.\;$  According to Arezzo-Tian \cite{Are-Tian03}
(c.f. \cite{song08} also), we can solve a homogenous complex equation starting
from nearby fibre:
\[
(\Omega+\sqrt{-1}\partial\overline{\partial}\Phi)^{n+1} = 0.
\]
Here $\Omega$ is the K\"ahler form in the total space $\M$ which is compatible to the
associated $C^{*}$ action of the test configuration.   We might scale in $\delta$ direction 
if necessary so that we solve this HCMA equation in $\M \rightarrow \triangle\;$ with $C^{\infty}$ regularities.  For $s\in [0,1]$, denote $\omega_s=(\Omega+\sqrt{-1}\partial\overline{\partial}\Phi)|_{M_s}$. The
geodesic ray $\phi_t$  starts from $\omega=\omega_1$, and it translates to the following equation in the space of
 complex structures compatible with a fixed symplectic form: $$\frac{dJ}{dt}=-J\mathcal DH,$$ where $H=\dot{\phi}_0$
 is a fixed function,  $J_t=f_t^*J|_{M_1}$, and $\dot{f}_t(x)=-\frac{1}{2}\nabla_t\dot{\phi}_t$. In this language,
  it is easy to see that the K-energy is convex along the ray:
 $$\frac{d^2 \bE}{dt^2}=\int_M{|\mathcal DH|^2\frac{\omega^n}{n!}}=\int_M{|\dot{J}|^2}\frac{\omega^n}{n!}.$$
 By definition,  the limit of the first derivatives of the K energy
 is precisely the $\yen$ invariant of this geodesic ray. According to \cite{chentang07}, this $\yen$ invariant
 agree with the Futaki invariant of the central fibre. By our assumption, the Futaki invariant of the central fibre vanishes.  Thus,
 \[
 \displaystyle \lim_{t \rightarrow \infty}\; {{d\; \bE}\over {d\, t}} = 0.
 \]
 Therefore, the K-energy will have a lower bound along the geodesic ray if we could prove the right hand side
  decays exponentially fast.
Now
the geodesic ray corresponds to
a foliation by holomorphic discs in $\M--$the total space of the underlying simple test configuration.
This foliation defines diffeomorphisms $f_s: M_1\rightarrow M_s(s\in [0,1])$, which satisfy $f_s^*\omega_s=\omega_1$.
 Then,  $\tilde{J}(s)=f_s^*(J|_{M_s})$ is a smooth family of complex structures compatible with $\omega_1$. Not surprisingly,
 this should coincide with the previous family $J_t$  up to a scaling of time parameter, i.e. $J_t=\tilde{J}_{e^{-t}}$.
 So $\dot{J}_t=-e^{-t}\dot{\tilde{J}}_{e^{-t}}$, and $|\dot{J}|^2\leq Ce^{-2t}$ for some constant $C>0$, which in turn
 proves the proposition.\\
 
 Note that the derivative of Calabi energy alone geodesic ray as
 \[
 \begin{array}{lcl} |{d\over {d\, t}} Ca(\omega_{\phi(t)})| & = & |\int_{M}\; H_{,\alpha\beta} \cdot R^{,\alpha \beta} \omega_{\phi(t)}^{n} |\\ & \leq & \sqrt{\int_{M} |D H|^{2} \omega_{\phi(t)}^{n}} \sqrt{\int_{M}\;|D R(\phi(t)|^{2} \omega_{\phi(t)}^{n}}\\
 & \leq & C e^{- t}
 \end{array}
 \] for some uniform constant $C.\;$   This suggests that the limit of the Calabi energy exists along any smooth geodesic ray associated to the simple test configuration. Set
 \[
 A = \displaystyle \lim_{t \rightarrow \infty}\; Ca(\omega_{\phi(t)}).
 \]
 Then, the difference $|Ca(\omega_{\phi(t)}) - A|$ will converges to $0$  exponentially fast. \\
 
 Smooth geodesic rays are constructed  \cite{Are-Tian03}  from nearby fibers.  From the way of construction,  it is hard to control where the limit of these geodesic ray is.  A new observation in \cite{song08} is that:
S. Sun can  specify which limit metric in the central fibre a priori;  then he showed that there always exists  a  smooth geodesic ray \footnote{He might have to chose a fibre much closer to the central fibre than what described in
 \cite{Are-Tian03}.} whose limit is this preferred metric in the central fibre.   If the central fibre has cscK metric,  it is convenient for us to choose a geodesic ray whose limit is this cscK metric.  For such a ray, we know that the limit of the Calabi energy must be
 $0.\;$ It follows that the Calabi energy  decays exponentially. \end{proof}

\section{Proof of Theorem 1.5}
 \subsection{Notations and set up}
 This subsection is a technical preparation of next subsection. \\

For any two K\"ahler potentials $\phi_0, \phi_1\in \cal H,$ we
want to use the almost smooth solution to approximate the $C^{1,1}$
geodesic between $\phi_0$ and $\phi_1.\;$ 
Let us setup some notations first.  Let $\Sigma^{(\infty)} = (-\infty,\infty) \times [0,1] \subset  {\bf R}^2$ denote the infinitely long strip.  For any integer $l, $ let $\Sigma^{(l)} $ be a long ``oval shape" disc such
that
$\Sigma^{(l)}$ is the union of $[-l,l]\times [0,1]$ with a half circle centered at $(-l, {1\over 2})$ with radius
${1\over 2}$ at the left, and a half circle centered at $(l, {1\over 2})$ with radius
${1\over 2}$ at the right.   Note that we want to smooth out the corner at the four corner points $\{\pm l\} \times \{0,1\}\;$ so that $\Sigma^{(l)}$ is a smooth domain\footnote{The author wish to stress that we do this ``smoothing" once for all: Namely, we smooth $\Sigma^{(1)}$ first. For any $l > 1$,  we may construct  $\Sigma^{(l)} (l\geq 1)$  by replacing the central line segment $\{0\}\times [0,1]$ in $\Sigma^{(1)}$ by  a cylinder $[-l+1, l-1] \times [0,1].\;$}. By construction, $\{\Sigma^{(l)}\}$ is a sequence of long ovals  in
this infinite strip $\Sigma^{(\infty)}\;$ where $\Sigma^{(0)}$ is a disc of radius ${1\over 2}$ centered at
$(0,{1\over 2}).\;$\\

Let $\psi$ be a convex family of K\"ahler potentials in $\Sigma^{(\infty)}$ given by
\begin{equation}
\psi(s,t, \cdot) =\bar \phi(s,t,\cdot) \qquad\qquad \forall (s,t) \in \Sigma^{(\infty)}.
\label{eq:stripeboundaryvalue}
\end{equation}
Here $\bar \phi(s,t,\cdot)$ can be any convex path connecting $\phi_0, \phi_1.\;$ For instance,
we may set
\[
\bar \phi(s,t, \cdot) = (1-t) \phi_0 + t \phi_1 - K t(1-t)
\]
where $K$ is a large enough constant. Here $K$ must depend on $\phi_0, \phi_1$ to ensure this family of potentials is convex.  In this subsection,  we assume that our boundary map $\psi$ is independent of $s$ variable.\\

Consider
Dirichlet problem for the HCMA equation \ref{eq:hcma0} on the
long oval shape domain $\Sigma^{(l)} $ with boundary
value
\[
\phi\mid_{\p \Sigma^{(l)}\times M} =  \psi\mid_{\p \Sigma^{(l)}\times M}.
\]
As in \cite{chen991}, we want to solve this via approximation method.
For any $\epsilon > 0$, consider the Drichelet problem:

\begin{equation}
(\pi_2^*\omega + \p \bar \p \phi)^{n+1} =  \epsilon  \cdot (\pi_2^*\omega + \p \bar \p \psi)^{n+1},
\qquad \forall\; (s,t,\cdot) \in  \Sigma^{(l)} \times M
\label{eq:leequation11}
\end{equation}
with fixed boundary data
\begin{equation}
\phi\mid_{\p \Sigma^{(l)} \times M} =  \psi\mid_{\p \Sigma^{(l)}\times M}.
\label{eq:leequation01}
\end{equation}
Denote the solution to this Dirichlet problem as $\phi^{(l,\epsilon)}$ for any $l\geq 1$ and $\epsilon\in (0,1).\;$
For each fixed $l$, one can prove as in [9] that solutions $\phi^{(l,\epsilon)}$ to eq.(4.2) (4.3)  have $C^{1,1}$ upper bounds independent of $\epsilon$. 

\begin{theo} \cite{chen05} For every $l$ fixed, there is a $C^{1,1}$ solution $\phi^{(l)}$ to the
equation (4.2)(4.3)  with $\epsilon =0.\;$ More importantly,
this upper bound on $ |\p \b \p \phi^{(l,\epsilon)}|$  is independent of $l>1\;$ and $\epsilon \in(0, 1].\;$
\end{theo}
\begin{theo} \cite{chen05} There exists a sequence of almost smooth solutions $\{\phi^{(l,\epsilon)}, l\in \NN, \epsilon \in (0,1)\}$ such that the upper bound of $ |\p \b \p \phi^{(l,\epsilon)}|$  is uniform and independent of $l>1\;$ and $\epsilon \in(0, 1].\;$ More importantly, when $(l,\epsilon) \rightarrow (\infty, 0)$, this sequence of almost smooth solutions converges to the unique $C^{1,1}$ geodesic connecting
the two end ``points'' ($\varphi_{0}, \varphi_{1}$).
\end{theo}

 Both theorems are crucial in the proof of Theorem 1.5 below. However, we will omit proof here.  In Section 6, we will give proof to some more general theorems in Section 5.
 \subsection{Proof of Theorem 1.5}
 The proof of Theorem 1.5 largely follows from the proof of Lemma 3.1. However, we have not proved that the space of K\"ahler potential is $C^\infty$ connected by geodesic segments yet. Thus, the proof
of Theorem 1.5 is more complicated since we need to overcome the lack of sufficient regularity.
\begin{proof}  Suppose $\rho:[0,\infty)\rightarrow\infty$ is a
geodesic ray parametrized by arc length $s$ such that
\[
  \displaystyle \lim_{s\rightarrow \infty}(d\,\bE, {{\p\rho}\over
  {\p s}})_{\rho(s)} < 0
\]
with
\[
\|{{\p \rho}\over {\p s}}\|_{\omega_{\rho(s)}}  = 1, \qquad \forall s \in [0, \infty).
\]
For any K\"ahler potential $\varphi_0\in \cH$, consider the unique
$C^{1,1}$ geodesic segment connecting $\varphi_0$ to $\rho(l).\;$  Suppose the length
of this geodesic segment is $\tau(l) (\forall l > 0).\;$  We  denote 
this geodesic segment as $ \Psi_l: [0, \tau(l)] \rightarrow \cH.\;$  For any $l$ large enough, we have
\[\begin{array}{lcl} &&
 (d\,\bE, {{\p\Psi_l}\over {\p
t}}\mid_{t=\tau(l)})_{\rho(l)} \\ & = & (d\,\bE, {{\p\Psi_l}\over {\p
s}}\mid_{t=\tau(l)} - {{\p \rho}\over {\p s}}\mid_{s=l} )_{\rho(l)} + (d\,\bE,{{\p \rho}\over {\p s}}\mid_{s=l}
)_{\rho(l)}\\ & \leq & \left(\int_{M}\; (R(\rho(l)) - \underline{R})^2
\omega_{\rho}^n\right)^{1\over 2} \cdot \left(  \int_M\;( {{\p\rho}\over {\p
s}}\mid_{s=l} -  {{\p\Psi_l}\over {\p
t}}\mid_{t=\tau(l)} )^2 \omega_{\rho(l)}^n\right)^{1\over 2}  + (d\,\bE,{{\p \rho}\over {\p s}}\mid_{s=l})_{\rho(l)}\\
& \leq  & Ca(\omega_{\rho(l)}) ^{1\over 2} \cdot (2 - 2 ( {{\p\rho}\over {\p
s}}\mid_{s=l},\; {{\p \Psi_l}\over {\p t}}_{t=\tau(l)})_{\rho(l)})^{1\over 2}+ (d\,\bE,{{\p \rho}\over {\p
s}}\mid_{s=l})_{\rho(l)}.
\end{array}
\]
With $l$ large enough, we can essentially treat $\tau(l) = l.\; $
According to Calabi-Chen\cite{chen992},  the infinite dimensional space $\cH$ is a non-positively curved manifold in the sense of
Alexandrov.  Thus, the small  angle at $\rho(l)$ on this  long, thin geodesic hinge approaches $0$ as  $l \rightarrow \infty.\;$  Moreover, it is smaller than the small angle of the corresponding
long, thin hinge in Euclidean plane.  Thus,  we have
\[
0  \leq l^2 \cdot \left( 1 - ( {{\p\rho}\over {\p
s}}\mid_{s= l},\; {{\p \Psi_l}\over {\p t}}_{t=\tau(l)} )_{\rho(l)} \right)  \leq 100 \;d(\varphi_{0},\rho(0))^{2}.
\]
Here $d(\varphi,\psi)$ denote the geodesic distance between two K\"ahler potentials $\varphi,\psi \in \cal H.\;$   Thus, we have
\begin{equation}
 (d\,\bE, {{\p\Psi_l}\over {\p
t}}\mid_{t=\tau(l)})_{\rho(l)} \leq { {10 \;\sqrt{Ca(\omega_{\rho(l)})} \cdot d(\varphi_{0}, \rho(0))} \over l}  + (d\,\bE,{{\p \rho}\over {\p
s}}\mid_{s=l})_{\rho(l)}.
\end{equation}
Now, we need to compare the K energy at $\varphi_{0}$ with $\rho(l).\;$ For convenience, let
us parametrize the geodesic segment between these two potentials in $[0,1].\;$ Then, the length
of ${{\p \Psi_{l}}\over {\p t}} $ is $\tau(l).\;$ \\

Following \cite{chen05}, we set  ${\kappa}:(-\infty,\infty) \rightarrow {\bf R}$ be a smooth
non-negative function such that ${\kappa} \equiv 1 $ on $[ -
{1\over 2},{1\over 2}]$ and vanishes outside  of $[-{3\over
4},{3\over 4}].\;$
Set\[ {\kappa}^{(l)} (s) = {1\over v} {\kappa}({s\over l}), \qquad {\rm where} \;v = \int_{-\infty}^\infty\; {\kappa}(s)\;d\,s. \]
 For any $m < l$, set
\[
  \bE^{(ml)}(t)  =  \displaystyle \int_{-\infty}^\infty\; k^{(m)}(s) \; \bE^{(l)}(s,t) \;d\,s, \qquad \forall m \leq l \in \NN
\]
and
\[
E^{(l)}(s,t) = \bE(\omega_{\phi^{(l,\epsilon)}}), \qquad \forall \;\;(s,t) \in \Sigma^{(l)}.
\]
For simplicity of presentation, we omit $\epsilon$ from our notations.
Then,
\[
  \bE^{(ml)}(0) =\bE(\varphi_0),  \; \bE^{(ml)}(1) = \bE(\rho(l)).
\]
Set
\[
 f^{(ml)}(t) =  {{d\, \bE^{(ml)}}\over {d\,t}}(t), \qquad \forall\; t\in [0,1].
\]
Following the same calculation in \cite{chen05}, for any $ 0\leq t_1 < t_2 \leq 1$
\begin{eqnarray}
f^{(ml)} (t_2) - f^{(ml)} (t_1) & \geq & -{1\over m^2} {1\over v}
\displaystyle \int_{t_1}^{t_2}\; \displaystyle \int_{-\infty}^\infty {{d^2
\kappa}\over {d\,s^2}}\mid_{s\over m}\;
\bE^{(l)}(s,t)\;d\,s\;d\,t\;\nonumber
\\
&\geq & - \displaystyle \int_{t_1}^{t_2}\;{1\over m^2} {1\over v}\;\mid
\displaystyle \int_{-\infty}^\infty {{d^2 {\kappa}}\over
{d\,s^2}}\mid_{s\over m}\; \bE^{(l)}(s,t)\;d\,s\mid
\;d\,t\nonumber
\\&\geq &- \displaystyle \int_{t_1}^{t_2}\;  {C\over
m} \;d\,t = -{C\over {2 m}}. \label{eq:firstderivativeofKenergy1}
\end{eqnarray}
Set $t_2 =1 $ and $t_1 = t \in [0,1].\;$ We have
\[
f^{(ml)}(t) \leq {C\over {2 m}} + f^{(ml)}(1).
\]
Therefore, 
\[
\begin{array}{lcl}  \bE(\rho(l)) - \bE(\varphi_0) & = &  \bE^{(ml)}(1) - \bE^{(ml)}(0)\\
& = & \displaystyle \int_0^1\;  {{d\, \bE^{(ml)}}\over {d\,t}}(t)  \; d\,t =  \displaystyle \int_0^1\; f^{(ml)}(t) \; d\,t\\
& \leq  &  \displaystyle \int_0^1\;\left( f^{(ml)}(1) + {C\over {2m}}\right)\; d\,t\\
& = &  \displaystyle \int_{-\infty}^\infty\; {\kappa}^{(m)}(s)
(d\,\bE, {{\p\Psi_l}\over {\p
s}}\mid_{s=\tau(l)})_{\rho(l)}  \cdot \tau(l) \;d\,s  + {C\over {2m}}\\
& = &  \displaystyle \int_{-\infty}^\infty\; {\kappa}^{(m)}(s) \left( { {10 \;\sqrt{Ca(\omega_{\rho(l)})} \cdot d(\varphi_{0}, \rho(0))} \over l}  \cdot \tau(l) + (d\,\bE,{{\p \rho}\over {\p
t}}\mid_{t=l})_{\rho(l)} \cdot \tau(l) \right)\;d\;s \\
 && \qquad +  {C\over {2m}}\\
& = &    10 \;\sqrt{Ca(\omega_{\rho(l)})} \cdot d(\varphi_{0}, \rho(0))\; {\tau(l) \over l}  + (d\,\bE,{{\p \rho}\over {\p
s}}\mid_{s=l})_{\rho(l)} \cdot \tau(l) +  {C\over {2m}}.
\end{array}
\]
Since the geodesic ray is destabilizing, so we know that
\[
(d\,\bE,{{\p \rho}\over {\p
s}}\mid_{s=l})_{\rho(l)} \leq 0.
\]
Thus, we have
\[
 \bE(\rho(l)) - \bE(\varphi_0) \leq 10 \;\sqrt{Ca(\omega_{\rho(l)})} \cdot d(\varphi_{0}, \rho(0))\; {\tau(l) \over l}  + {C\over {2m}}
\]
As before, let $l \rightarrow \infty$, we have

\[
E(\rho(l)) - E(\varphi_{0}) \leq {C\over {2m}}
\]
Let  $m \rightarrow \infty$, then we have
\[
E(\varphi_{0}) \geq \lim_{l\rightarrow \infty}\;E(\rho(l)) \geq 0.
\]
Thus the K energy is bounded from below in the entire K\"ahler class since $\varphi_{0}$ is an arbitrary
chosen K\"ahler potentials. 

\end{proof}

\section{Proof of Theorem 1.1}
\subsection{Setup and main results}

Let us setup some notation first. The setup is a slight modification of the setup in Section 4.1.   \\

Let $\Sigma^{(\infty, m)} = (-\infty,\infty) \times [0,m] \subset  \RR^2$ denote a sequence of infinitely long strips.  Fix $m \in \NN$, for any integer $l, $ let $\Sigma^{(l, m)} $ be a long ``oval shape" disc such
that
$\Sigma^{(l,m)}$ is the union of $[-l,l]\times [0,m]$ with a half circle centered at $(-l, {m\over 2})$ with radius
${m \over 2}$ at the left, and a half circle centered at $(l, {m \over 2})$ with radius
${m \over 2}$ at the right.   Note that we want to smooth out the corner at the four corner points $\{\pm l\} \times \{0,m\}\;$ so that $\Sigma^{(l,m)}$ is a smooth domain\footnote{The author wish to stress that we do this ``smoothing" once for all: Namely, we smooth $\Sigma^{(1,m)}$ first. For any $l > 1$,  we may construct  $\Sigma^{(l,m)} (l\geq 1)$  by replacing the central line segment $\{0\}\times [0,m]$ in $\Sigma^{(1,m)}$ by  a cylinder $[-l+1, l-1] \times [0,m].\;$}. By construction, $\{\Sigma^{(l, m)}\}$ is a sequence of long ovals  in
this infinite strip $\Sigma^{(\infty,m)}\;$ where $\Sigma^{(0,m)}$ is a disc of radius ${m\over 2}$ centered at
$(0,{m\over 2}).\;$\\

Given $\varphi_{0}$ and a degenerating  geodesic ray $\rho:  [0,\infty)  \rightarrow \cH.\;$ By assumption, this geodesic ray is tamed by a bounded ambient metric 
\[
h = \pi_{2}^{*} \omega + \sqrt{-1}\;\bar \partial \partial \bar \rho
\]
in $[0,\infty) \times S^{1}\times M.\;$  Let $\bar \phi$ be any convex path connecting $\varphi_{0}$ and $\bar \rho(0).\;$  Putting this segment together with the initial ray $\bar\rho([0,\infty))$,
adjusting it by  first smoothing out the corner and then adding a function of $\tau$ only,  we can obtain a one parameter convex family of K\"ahler potentials. We denote this family as $\bar \phi: [0,\infty) \rightarrow \cH.\;$ Then, $\bar \phi$ is essentially same as the initial geodesic ray $\bar \rho.\;$  By abusing notations, we may set
\[
\bar \phi = \bar \rho.
\] 
 Note that the K energy of $\omega_{\bar \rho (\tau,\cdot)}  (\tau \in [0,\infty))$ is uniformly bounded from above and below. This is
 an important fact which will be used crucially later.\\
 
As in \cite{chen991} \cite{chen05}, we want to solve the disc version geodesic problem via approximation method.
For any $\epsilon > 0$, consider the Drichelet problem:

\begin{equation}
(\pi_2^*\omega + \p \bar \p \phi)^{n+1} =  \epsilon  \cdot (\pi_2^*\omega + \p \bar \p \bar\rho)^{n+1},
\qquad \forall\; (s,t,\cdot) \in  \Sigma^{(l,m)} \times M
\label{eq:leequation11}
\end{equation}
with fixed boundary data
\begin{equation}
\phi\mid_{\p \Sigma^{(l,m)} \times M} = \bar \rho \mid_{\p \Sigma^{(l,m)}\times M}.
\label{eq:leequation01}
\end{equation}
Denote the solution to this Dirichlet problem as $\phi^{(l,m,\epsilon)}$ for any $l,m \gg 1$ and $\epsilon\in (0,1).\;$ Following \cite{chentian005}, set
\[
\phi^{(l,m)} = \displaystyle \lim_{\epsilon \rightarrow 0}\; \phi^{(l,m,\epsilon)}
\]
where $\{\phi^{(l,m,\epsilon)}\}$  is a sequence of almost smooth solution in $\Sigma^{(l,m)}.\;$  Set
\[
\phi^{(m)} = \displaystyle \lim_{l \rightarrow \infty}\; \phi^{(l,m)}
\]
 is a  $C^{1,1}$ geodesic segment between $\varphi_0$ and $\rho(m).\;$
Set
\[
\phi = \displaystyle \lim_{m \rightarrow \infty}\; \phi^{(m)}: [0,\infty) \rightarrow \cH.
\]
Then, $\phi$ is a geodesic ray initiating from $\varphi_0$ but parallel to $\rho.\;$ Thus,
\begin{equation}
\displaystyle \sup_{m \in \NN} \displaystyle \max_{[0,m)\times M}  \left( \|\phi^{(m)} - \phi\|_{C^{0}}  + \| \rho- \bar \rho\|_{C^{0}} + \|\phi- \rho\|_{C^{0}}\right) \leq C.
\end{equation}

For each fixed $l,m$, one can prove as in [9] that the almost smooth solutions $\phi^{(l,m,\epsilon)}$ to eq.(\ref{eq:leequation11})
 (5.3) have a uniform $C^{1,1}$  upper bounds (independent of $\epsilon$).  The main technical theorem in this section is to show 

\begin{theo} For any geodesic ray which is tamed by a bounded ambient geometry,  for every $(l,m,\epsilon) \in (1,\infty) \times (1, \infty) \times (0,1),$  there exists a sequence of almost smooth solution
$\phi^{(l,m,\epsilon)}(t, s, \cdot)$ 
\[
   \begin{array}{rcl} (\pi_{2}^{*} \omega + \sqrt{-1} \partial \b \partial \phi^{(l,m,\epsilon)})^{n+1} & = & \epsilon ( \pi_{2}^{*} \omega + \sqrt{-1} \partial \b \partial \bar \rho)^{n+1} \\
   \phi^{(l,m,\epsilon)} \mid_{\p \Sigma^{(l, m)}} & = & \bar \rho,\end{array}
\] 
such that every $m$ fixed, let $(l, \epsilon) \rightarrow (\infty, 0)$, the sequence of K\"ahler potential
approximate to a $C^{1,1}$ geodesic segment $ \phi^{(m)}.\;$ When $m \rightarrow \infty$, this family of geodesic segments converges to the relative $C^{1,1}$ geodesic ray $\phi$ initiated from $\varphi_{0}$ in the direction of the initial geodesic ray.  More importantly,
\[ \; |\p \b \p \left(\phi^{(l,m,\epsilon)}(t,s,\cdots) - \bar \rho\right)|_{h} \leq C\]  
where $C,$
are independent of $l, m >1\;$ and $\epsilon \in(0, 1].\;$ \end{theo}

We delay the proof of this theorem.  As an application, we prove

\begin{theo}  If $\varphi_{0}$ is a cscK metric,  then $\{\phi^{(l,m,\epsilon)} \}$ converges
to a relative geodesic ray $\phi$ initiated from $\varphi_{0}$ such that every K\"ahler metric $\omega_\phi$  in this ray is
a $C^{1,1}$ K-energy minimizer.  Moreover, 
the convergence of ${\omega_{\phi^{(l,m,\epsilon)}}^n \over \omega^n}$ to ${\omega_\phi^n \over \omega^n}$  is strong in $L^{2}(M, \omega)$ norm.                                                            
\end{theo}

\begin{theo} Any relative $C^{1,1}$ geodesic ray consists of K-energy minimizers only
must be also smooth if  this geodesic ray is a limit of a sequence of almost smooth solutions to the disc version of geodesic equation. 
\end{theo}

Proof of Theorem 5.1 and 5.2  will be delayed to later subsections.  The proof of Theorem 5.3 is very similar to the corresponding results in \cite{chentian005} (after Theorem 5.2) is proved).  We will omit its proof here. Now we give a proof of Theorem 1.1.  
\begin{proof} Let $\varphi_0$ be a cscK metric and $\rho: [0,\infty)\rightarrow \cH$ be a degenerating geodesic ray tamed by a bounded ambient geometry.   Let $\phi^{(l,m,\epsilon)}$ be a family
of almost smooth solutions in $\Sigma^{(l,m)}$ with appropriate boundary data. Then,
\[
\phi = \displaystyle \lim_{m \rightarrow \infty} \; \displaystyle \lim_{l \rightarrow \infty} \;\displaystyle \lim_{\epsilon \rightarrow 0} \;\phi^{(l,m,\epsilon)}
\]
is a relative $C^{1,1}$ geodesic ray initiating from $\varphi_0$ but parallel to the initial geodesic ray
$\rho.\;$ According to Theorem 5.2, we know that for any $t\in [0,\infty)$, we have
\[
\bE (\omega_{\phi(t,\cdot)}) =  \displaystyle \inf_{\psi \in [\omega] } \;\bE(\omega_\psi)
\]
is a K-energy minimizer.  Since $\phi(t,\cdot)$ is  the $C^{1,1}$ limit of a sequence of almost smooth solutions, Theorem 5.3 implies that $\omega_{\phi(t,\cdot)}$ is smooth K\"ahler form for any $t\in [0,\infty).\;$ Thus, $\omega_\phi$ has constant scalar curvature.  It follows that $\phi: [0,\infty) \rightarrow \cH$
is a smooth geodesic ray consisting of cscK metrics only.   In particular, this means that the $\bE(\omega_{\phi(t,\cdot)} )$ is a constant function on $t.\;$By a direct calculation, we have
\[
0 = {{d^2 \bE (\omega_{\phi(t,\cdot)} )}\over {\p t^2}} = \displaystyle \int_M \; 
\| \bar \partial \nabla^{(1,0)} {{\partial \phi(t,\cdot)} \over {\p t}}\|^2_{\omega_{\phi(t,\cdot)}} \; \omega_{\phi(t,\cdot)}^n.
\]
Then,
\[
\bar \partial \nabla^{(1,0)} {{\partial \phi(t,\cdot)} \over {\p t}} = 0.
\]
It follows that
$\nabla^{(1,0)} {{\partial \phi(t,\cdot)} \over {\p t}}$ is a holomorphic vector field.  Therefore, $\phi:[0,\infty)\rightarrow \cH$ represents a holomorphic line in $\cH.\;$ The initial geodesic ray $\rho:[0,\infty) \rightarrow \cH$ is parallel to this holomorphic line by definition of $\phi.\;$ Our theorem is then proved.
\end{proof}
\subsection{Proof of Theorem 5.2}
 We need to prove a few Lemmas  first.  We follow notations from previous subsection.
\begin{lem}  Let $\rho: [0,\infty)\rightarrow  \cal H$ be a degenerated geodesic ray tamed by a bounded ambient geometry.  If the infimum of the Calabi energy along this ray is $0$\footnote{alternatively, we may assume  that there exists a cscK metric and a geodesic ray with bounded K-energy .},  then $\bE(\phi^{(l,m,\epsilon)})$ are uniformly bounded among this sequence of almost smooth solutions
in any  domain $\Sigma^{(l,m)}.\;$
\end{lem}
\begin{proof}

Since the K-energy functional is sub-harmonic over solution to disc version of geodesic equation, we have
\[
{{\p^{2} \bE(\phi^{(l,m,\epsilon)})}\over {\p z \p \bar z}} \geq 0.
\]
At the boundary, we know
\[
\bE(\phi^{(l,m,\epsilon)}(s,t,\cdot ) )= \bE(\bar \rho(t,\cdot)\mid_{\p \Sigma^{(l,m)}}) \leq C.
\]
By the maximal principle,  we know that $ \bE(\phi^{(l,m,\epsilon)}) $ is uniformly bounded from above. Here the upper bound is independent of $l,m, \epsilon.\;$\\

 If we assume that there exists a cscK metric, then the K-energy functional in $[\omega]$ has a uniform lower bound.  Alternatively, if we assume that $\rho: [0,\infty)\rightarrow \cal H$ is a degenerating geodesic ray where the Calabi energy approaches to $0$ in this ray, Theorem 1.5 also implies that the K-energy  is bounded from below in $(M, [\omega]).\;$  In either case, we can apply these a priori estimate to prove that  the K-energy functional $\bE(\phi^{(l,m,\epsilon)})$ is uniformly bounded from below in $(M, [\omega]). $ Our lemma is then proved.

  \end{proof} 
 
 \begin{lem} If there exists a cscK metric and a  degenerating geodesic ray where the Calabi energy converges to $0$ in this ray, then the first derivative of the K-energy functional along the two ``long side" is very small. More specifically, we have
\[
{{\p \bE(\phi^{(l,m,\epsilon)}) }\over {\p t}} \mid_{t =0} = 0, \]
 and for any  $s \in (- {{3 l}\over 4},  {{3 l}\over 4}),\; $  we have\[ 
0\leq {{\p \bE(\phi^{(l,m,\epsilon)}) }\over {\p t}} \mid_{t =m} \leq C\cdot  \sqrt {ca(\omega_{\rho(m)})} \cdot {{d(\varphi_{0},\rho_{0})}\over m}
\]
for some fixed positive constant $C.\;$
\end{lem}
 
\begin{proof}   Let $\rho: [0,\infty) \rightarrow \cH$ be a degenerating geodesic ray where the infimum of the Calabi energy in this ray is $0$.  Without loss of generalities, set $\rho_{m} = \rho(m) $
such that $\displaystyle \liminf_{m\rightarrow \infty} Ca(\rho_{m}) = 0$ and
\[
 {{d\, \bE_{\omega}}\over {d\, \tau}}\mid_{\tau = m} = - (R(\omega_{\rho_{m}}) - \bar R,  {{\p \rho}\over {\p \tau}}\mid_{\tau = m})_{\rho_{m}} \leq 0.
\]
Set
\[
u^{(m)} = {{\p \phi^{(m)}}\over {\p t}}\mid_{t=0},\qquad v^{(m)} = {{\p \phi^{(m)}}\over {\p t}}\mid_{t=m}
\]
Set 
\[
w^{(m)} = {{\p \rho}\over {\p \tau}}\mid_{\tau=m}.
\]

Following Lemma 4.10  in \cite{chen05},  we have
\[
- (R(\phi^{(m)}) - \bar R, u^{(m)})_{\varphi_0}  \leq  - (R(\phi^{(m)}) - \bar R, v^{(m)})_{\rho_m}.  
\]

Following the fact this is a degenerating geodesic ray, we have 
\[
- (R(\rho_m) - \bar R, w^{(m)})_{\rho_m}  \leq 0 \qquad {\rm and}\; \int_{M}\; (R(\rho_{m}) -\bar R)^{2}\;\omega_{\rho_{m}}^{n} \rightarrow 0.
\]
If $\varphi_{0}$ is cscK metric, we have
\[
 - (R(\phi^{(m)}(0)) - \bar R, u^{(m)})_{\varphi_0}  = - (R(\varphi_{0}) - \bar R, u^{(m)})_{\varphi_0}  = 0.
\]
Thus, 
\[
\begin{array}{lcl}
0 &\leq & - (R(\phi^{(m)}) - \bar R, v^{(m)})_{\rho_m}  \\
& = & - (R(\phi^{(m)}) - \bar R, w^{(m)})_{\rho_m}  - (R(\phi^{(m)}) - \bar R, v^{(m)} -w^{(m)})_{\rho_m}\\
& \leq  &    - (R(\phi^{(m)}) - \bar R, v^{(m)} -w^{(m)})_{\rho_m} \leq \sqrt {ca(\omega_{\rho(m)})} \cdot |v^{(m)}-w^{(m)}|_{\rho_m}
\\
& \leq & C \cdot  \sqrt {ca(\omega_{\rho(m)})} \cdot {{d(\varphi_{0},\rho_{0})}\over m}
\end{array}\]
for some positive constant $C > 0.\;$\\

Since $\phi^{(l,m,\epsilon)}$ converges to $\phi^{(m)}$ uniformly in $C^{1,\alpha}$ norm as $l \rightarrow \infty$ and $\epsilon \rightarrow 0.\;$ Our lemma follows
from the last inequality above.
\end{proof}

\begin{lem} \label{appl:kenergyisapproharmonic}As $l,m$ are sufficiently large we may control the $L^1$ measure of  $\Delta_{s,t} \bE_\omega^{(l,m,\epsilon)}$ in $\Sigma^{(l,m)}$ as
\[
\displaystyle \int_{\Sigma^{(l,m)}}\; \triangle_{s,t}\;\bE_{\omega}(\phi^{(l,m,\epsilon)}(s,t,\cdot)) \;d\,s\,d\,t \leq 
d(\varphi_{0}, \rho_{0}) \cdot {l\over m} \cdot \sqrt{Ca(\omega_{\rho_{m}})} + C \cdot {m\over l}
\]
for some uniform constant $C.\;$
\end{lem}
When there is no confusion, we will drop the superscript $l, \epsilon$. Set
\[
\bE^{(l,m)} = \bE(\phi^{(l,m,\epsilon)}).
\]
Here we have suppressed the dependency on $\epsilon.\;$
\begin{proof}
Let $\xi:(-\infty,\infty) \rightarrow \RR$ be a smooth
non-negative cut-off function such that $\xi \equiv 1 $ on $[ -
{1\over 2},{1\over 2}]$ and vanishes outside $[-{3\over 4},{3\over
4}].\;$
\[
\begin{array}{lcl} &&\int_{t=0}^m \;\int_{s=-{l\over 2}}^{l\over
2} \; |\Delta_{s,t}\; \bE^{(l,m)}(s,t)|\;d\,s\;d\,t \\
& \leq &
\int_{t=0}^m \;\int_{s=-{l}}^{l} \xi({s\over l}) \Delta_{s,t}\;
\bE^{(l,m)}(s,t)\;d\,s\;d\,t\\
&=& \int_{s=-l}^l {{\p \bE^{(l,m)}(s,t)}\over {\p t}}\mid_0^m\;
\xi({s\over l})\;d\,s - {1\over l} \int_{t=0}^m
\;\int_{s=-{l}}^{l}
\xi'({s\over l}) {{\p \bE^{(l,m)}(s,t)}\over {\p s}} \;d\,s\;d\,t\\
&=& l\cdot \sqrt {Ca(\omega_{\rho(m)})} \cdot {{d(\varphi_{0},\rho_{0})}\over m}+ {1\over l^2}  \int_{t=0}^m \;\int_{s=-{l}}^{l}
\xi''({s\over l}) \bE^{(l,m)}(s,t) \;d\,s\;d\,t \\
&\leq &l\cdot  \sqrt {Ca(\omega_{\rho(m)})} \cdot {{d(\varphi_{0},\rho_{0})}\over m}+  {1\over l^2}  \int_{t=0}^m \;\int_{s=-{l}}^{l}
|\xi''({s\over l})|\cdot |\bE^{(l,m)}(s,t)| \;d\,s\;d\,t\\
&  \leq & l\cdot  \sqrt {Ca(\omega_{\rho(m)})} \cdot {{d(\varphi_{0},\rho_{0})}\over m} + {1\over l^2}
\int_{t=0}^m \;\int_{s=-{l}}^{l} C \;d\,s\;d\,t\\
& = &l\cdot \sqrt {Ca(\omega_{\rho(m)})} \cdot {{d(\varphi_{0},\rho_{0})}\over m} + C \cdot {m \over l}.
\end{array}
\]
\end{proof}

An immediate Corollary is:
\begin{cor} Let $l = {m \over \left(ca(\omega_{\rho_{m}})\right)^{1\over 4}} $ and $m \rightarrow \infty$, then
\[
\displaystyle \int_{\Sigma^{({l\over 2},m)}}\;  \triangle_{s,t}\;\bE_{\omega}(\phi^{(l,m,\epsilon)}(s,t,\cdot)) \;d\,s\,d\,t \rightarrow 0.
\]

\end{cor}

Set
\[
A = \displaystyle \inf_{\varphi\in \cH}\; \bE(\varphi).
\]
Then, 
\begin{lem}  \label{appl:kenergyisalmostconstant}For any point $(s,t)$ in a fixed compact domain $\Omega$ in
$\Sigma^{(l,m)}$, except perhaps a set of measure $0$, we have \[
\displaystyle\lim_{l,m\rightarrow \infty} \;\bE^{(l,m)}(s,t) = \displaystyle
\lim_{l\rightarrow \infty}\; \bE(\phi^{(l)}(s,t)) = A.\; \]
\end{lem}

\begin{proof} Set $f^{(l,m)} = \Delta_{s,t} \; \bE^{(l,m)}(s,t)\geq 0.\;$ In $\Sigma_{({l\over 2},m)} \subset \Sigma_{(l,m)}$,
we have \[ \displaystyle \lim_{l,m\rightarrow \infty}\displaystyle
\int_{\Sigma_{({l\over 2}, m)}} \;f^{(l,m)} = 0.\]
 Next, we decompose
$\bE^{(l,m)}$ into two parts:
\[
\bE^{(l,m)} = u^{(l,m)} + v^{(l,m)}, \qquad \;{\rm in}\qquad \Sigma_{({l\over 2}, m)}
\]
such that
\[
\left\{ \begin{array}{ll} \Delta_{s,t} \; u^{(l,m)} \; = \; 0, & {\rm where}\;  u^{(l,m)}\mid_{\p \Sigma_{({l\over 2},m)}} =\bE_{\omega}^{(l,m)},\\
 \Delta_{s,t}\; v^{(l,m)} \; =\;f^{(l,m)} \geq 0, &\; {\rm where}\; v^{(l,m)}\mid_{\p \Sigma_{({l\over 2},m)}}=0.\end{array}
\right.
\]
It is clear that $v^{(l,m)}\leq 0.\;$   Since $\bE_{\omega}^{(l,m)}$ is uniformly bounded, then $u^{(l,m)}$ is a uniformly bounded
harmonic function in $\Sigma^{({l\over 2},m)}$ such that
\[
  u^{(l,m)}(s,0) =A, \qquad \forall \; s\in [-{l\over
2},{l\over 2}].
\]
Taking limit as $l,m\rightarrow \infty, $ then $\{u^{(l,m)}\}$ will converge locally smoothly to a bounded harmonic function in half plan such that it is constant $A$ along the line $t=0.\;$ Then, such a function must be constant globally. In other words,   in any compact
sub-domain $\Omega \subset \Sigma^{(l,m)}$ fixed,  we have $\displaystyle\lim_{l,m \rightarrow
\infty}\; u^{(l,m)} =A.\;$ Consequently,
\[\begin{array}{lcl} A &\leq &
\displaystyle\limsup_{l,m\rightarrow \infty} \;\bE^{(l,m)} \\
& = & \displaystyle\limsup_{l,m\rightarrow \infty}\;(u^{(l,m)} +
v^{(l,m)})
\\&\leq & \displaystyle\lim_{l,m\rightarrow \infty}\; u^{(l,m)} = A.
\end{array}\]
Therefore, for every point in $\Omega$ (fixed), we have
\[
\displaystyle\lim_{l,m\rightarrow\infty} \bE^{(l,m)} =
\displaystyle\limsup_{l,m\rightarrow\infty} \bE^{(l,m)}= A.
\]

\end{proof}
The leading term of the K-energy functional is 
\[
\int_M\; \log {{\omega_\phi^n}\over {\omega^n} } \cdot   {{\omega_\phi^n}\over {\omega}^n} \cdot \omega^{n}.
\]
The rest part of the K energy functional converges strongly under weakly $C^{1,1}$ topology. If
we set
\[
f = \log {{\omega_\phi^n}\over {\omega^n}}.
\]
Then, the convergence of  K energy with uniform $C^{1,1}$ bound on K\"ahler potential implies
that \[
\int_M\; \log {{\omega_{\phi^{(l,m,\epsilon)}}^n}\over {\omega^n} } \cdot   {{\omega_{\phi^{(l,m,\epsilon)}}^n}\over {\omega}^n} \cdot \omega^{n}
\]
converges to the corresponding limit.
The third part of Theorem 5.2 (strong convergence of volume form ration in $L^{2} (M, \omega)$ space) then follows from 

\begin{lem} \label{appl:volumestrongconvergence} Suppose that $\{f_m, m \in \NN\}$ is a sequence of positive, uniformly bounded functions such that $f_{m}$ weakly converge to $g$ in $L^{p}$ for some $p$ large enough.
Then, $\{f_{m}\}$ converges strongly to $g$ in $L^{2}(M)$ if 
\[
\lim_{m\rightarrow \infty}\; \displaystyle \int_{M}\; f_{m} \log f_{m} = \displaystyle \int_{M}\; g\;\log g. 
\]
\end{lem}
\begin{proof}  Without loss of generality, may assume
\[
\displaystyle \int_{M}\; f_{m} =  \displaystyle \int_{M} g  = 1.
\]
By our assumption,  we have
\begin{equation}
 \int_M\; f_l \log f_l \;\omega^n - \int_M \;g \log g
= o({1\over l}).  \label{eq:strongvolumeconvergence}
\end{equation}
 Define
$F(u) = u \log u.\;$ For any $l$ large enough and for any
$\epsilon>0$, set
$$F(t) = F(t f_l + (1-t) (g+\epsilon)) = F(a\,t +
b),$$
where
\[
  a =  f_l - g -\epsilon, \qquad {\rm and}\qquad b = g + \epsilon.
\]
Note that $a,b $ are both functions in $M.\;$ Clearly, we have
\[ |a| + |b| \leq C.\;\]
Note that
\[
F'(t) = a \log (a\,t + b) + a,
\]
and
\[
F''(t) = {{a^2}\over {a t + b}} \geq {a^2\over C}, \qquad \forall
\; t \in [0,1].
\]
Thus,
\[
\begin{array}{lcl} \int_M F'(0)\omega^n& =& \int_M\; (a \log b + a)\omega^n\\ & =  &
\int_M\; (f_l-g-\epsilon) \log (g +\epsilon) \omega^n   +
\int_M\; (f_l-g-\epsilon) \omega^n\\ & = & \int_M\;
(f_l-g-\epsilon) \log (g +\epsilon) \omega^n    -\epsilon  
\end{array}
\]
Taking the following double limits
\[
\begin{array}{lcl}
\displaystyle\;\lim_{\epsilon \rightarrow 0} \displaystyle\;
\lim_{l\rightarrow \infty} \displaystyle\; \int_M\; F'(0) \omega^n
& = & \displaystyle\;\lim_{\epsilon\rightarrow 0}
\displaystyle\;\lim_{l\rightarrow \infty} \left(\int_M\;
(f_l-g-\epsilon)\; \log (g+\epsilon) \;\omega^n -\epsilon
\right)\\ &&
 \displaystyle\; \lim_{\epsilon\rightarrow 0}\; \epsilon =0.
\end{array}
\]
It follows,
\[
\begin{array}{lcl} F(1)-F(0) & =& F(f_l) - F(g+\epsilon) \\
& = & \int_0^1\; F'(t) d\,t = F'(0) + \int_0^1 \int_0^t
\;F''(s)\;d\,s\,d\,t \\
& = & F'(0) +  \int_0^1 \int_0^t \;{{a^2}\over {a s + b}}
\;d\,s\,d\,t\\
&\geq & F'(0) +   \int_0^1 \int_0^t \;{{a^2}\over { C }}
\;d\,s\,d\,t =  F'(0)  + {a^2\over {2 C}}.
\end{array}.
\]
Integrating this over $M,\;$ we have,
\[\begin{array}{lcl} &&
\int_M\; (f_l-g-\epsilon)^2 \omega^n  =  \int_M \;a^2 \omega^n
\\ \qquad \qquad & \leq & 2 C \int_M\; (F(1) - F(0))\omega^n -  2 C \int_M F'(0)
\omega^n \\ \qquad \qquad & = & 2 C \left(\int_M f_l \log f_l \;
\omega^n - \int_M \;(g +\epsilon) \log(g+\epsilon) \omega^n\right)
-  2 C \int_M F'(0) \omega^n 
\\ &\leq & o(\epsilon + {1\over l}) -  C \int_M F'(0) \omega^n. \end{array}
\]
Consequently, we have
\[
\begin{array}{lcl} 2 \displaystyle \lim_{l\rightarrow 0} \int_M\;(f_l-g)^2
\omega^n 
& \leq & \displaystyle \lim_{\epsilon\rightarrow 0} \displaystyle
\lim_{l\rightarrow \infty} \displaystyle  o({1\over l} +\epsilon)
+ C \displaystyle\;\lim_{\epsilon \rightarrow 0} \displaystyle\;
\lim_{l\rightarrow \infty} \displaystyle\; \int_M\; F'(0) \;\omega^n\\
& = & 0.
\end{array}
\]
Therefore, $f_{m} $ converges
strongly to $g$ in
$L^2(M,\omega).\;$
\end{proof}
\subsection{ Proof of Theorem 5.1}
The proof of Theorem 5.1 is very similar to Theorem 3.9 in \cite{chen05}.   As in \cite{chen991}, the proof of the existence of $C^{1,1}$ geodesic segment also prove the existence of $C^{1,1} $ solution for disc version geodesic solution as well. However, Theorem  5.1 in this paper and theorem 3.9 in \cite{chen05}
is more like is a non-compact  analogue statement of  problems addressed  in \cite{chen991}. We want to discuss
more carefully here.  Nonetheless, the key estimate should be similar conceptually, especially when $m \ge 1$.\\

As in \cite{chen991}, let $\phi^{(m,\epsilon)}(t, \cdot)$ denote the $\epsilon$-approximated $S^{1}$ invariant  solution
for geodesic equation between $\varphi_0$ and $\rho_m.\;$   Since $\phi^{(m,\epsilon)}$ is independent
of $s$, we can view it as a solution in $\Sigma^{(\infty,m)}\times M.\;$  In other words, we have 
\begin{equation}
(\pi_2^*\omega + \p \bar \p \phi)^{n+1} =  \epsilon  \cdot (\pi_2^*\omega + \p \bar \p \bar \rho)^{n+1},
\qquad \forall\; (t,\cdot) \in  \Sigma^{(\infty,m)}\times M \label{eq:leequation21}
\end{equation}
with Dirichlet boundary data
\[
\phi^{(m,\epsilon)}(s, 0,\cdot) = \varphi_0, \;\; \phi^{(m,\epsilon)}(s, m,\cdot) = \rho_m.
\]
Here we  abuse notations by letting
\[
\phi^{(m,\epsilon)}(s,t,\cdot) = \phi^{(m,\epsilon)}(t, \cdot).
\]
Let $\phi^{(m)}$ denote the $C^{1,1}$ geodesic between $\varphi_0, \rho_m.\;$ By Maximum Principle, we
know that $\phi^{(m,\epsilon)}$ monotonically increases as $\epsilon $ decreases to $0.\;$ In particular,
we have
\begin{equation}
\bar \rho(s,t,\cdot) \leq \phi^{(m,\epsilon)}(s,t,\cdot) \leq \phi^{(m)}(t,\cdot), \qquad \forall (s,t) \in \Sigma^{(\infty,m)}.
\label{eq:elequation3} \end{equation}
In particular, there is an error term $o(\epsilon)$ such that $\displaystyle \lim_{\epsilon \rightarrow 0}\; o(\epsilon) = 0$ in any $C^{1,\alpha}$ norm and
\begin{equation}
\phi^{(m,\epsilon)}(s,t,\cdot) = \phi^{(m)}(t) + o(\epsilon). \label{eq:elequation41} 
\end{equation}

For any $l \in \NN$, choose $\delta_l$ such that
\[
\displaystyle \lim_{l \rightarrow \infty}\; (l-\delta_l)  = \infty.
\]

\begin{lem}   The sequence of almost smooth solution $\{\phi^{(l,m, \epsilon)}\}$
converges uniformly continuous to $\phi\;$ in $\Sigma^{(\delta_l, m)} \times M$ when $l, m \rightarrow \infty$ and $\epsilon \rightarrow 0.\;$  In particular, the same convergence result
holds for $\{\phi^{(l,m)}\}$ in any sub-domain $\Sigma^{(\delta_l,m)}\times M.\;$
\end{lem}

To prove this Lemma, we need to introduce  a sequence of harmonic functions
$h^{(l,m)}$ in $\Sigma^{(l,m)}$ such that its boundary value of $h^{(l,m)}$ in $\p \Sigma^{(l,m)}$ is
\[
h^{(l.m)} (s,t) =\left\{ \begin{array}{ll} 0 & |s| \leq l-1,\\
K & |s| \geq  l,\\
\in [0,K] & {\rm otherwise} \end{array} \right. 
\]
where $K$ is some large enough positive constant:
\[
K >  2 \displaystyle \sup_{m \in \NN}  \displaystyle \max_{[0,m] \times M} \; |\phi^{(m)} -\bar \rho| + 1.
\]
  The following Lemma is critical
\begin{lem} In $\Sigma^{(\delta_l, m)} \times M$, we have
\[
\displaystyle  \lim_{l,m\rightarrow \infty} \displaystyle \max_{\Sigma^{(\delta_l, m)} \times M} \;h^{(l,m)} = 0.
\]
\end{lem}
The proof, which we omit, is elementary.  Now we are ready to prove Lemma 5.10. \\

\begin{proof}
Note that for any $1\geq \epsilon > 0$ and $l$ fixed, we have
\[
\begin{array}{lcl}
\left(\pi_2^*\omega + \p \bar \p \phi^{(l,m,\epsilon)}\right)^{n+1} &  = &  \epsilon  \cdot (\pi_2^*\omega + \p \bar \p \bar \rho)^{n+1}\\
& = & (\pi_2^*\omega + \p \bar \p \phi^{(m,\epsilon)})^{n+1} 
\end{array}
\]
with boundary data
\[
\begin{array}{lcl}
 \phi^{(l,m,\epsilon)}\mid_{\p \Sigma^{(l,m)} \times M} &= & \bar \rho\mid_{\p \Sigma^{(l,m)}\times M}
 \\ & \leq & \phi^{(m,\epsilon)}\mid_{\p \Sigma^{(l,m)} \times M} 
 \end{array}
\]
By the maximum principle for Monge-Ampere type equation, we have
\[
\phi^{(l,m, \epsilon)} \leq \phi^{(m,\epsilon)}, \qquad \qquad \forall\;\; (s,t,\cdot) \in \Sigma^{(l,m)}\times M.
\]
Combining this with inequality (\ref{eq:elequation3}), we have
\[
\phi^{(l, m, \epsilon)} \leq \phi^{(m)}, \qquad \qquad \forall\;\; (s,t,\cdot) \in \Sigma^{(l,m)}\times M.
\]
It is easy to see that for any $(l, m, \epsilon)$ we have
\[
\bar \rho(s,t,\cdot) \leq \phi^{(l, m, \epsilon)} , \qquad \qquad \forall\;\; (s,t,\cdot) \in \Sigma^{(l,m)}\times M.
\]

Note that $h^{(l,m)}$ is a harmonic function on $\Sigma^{(l,m)}$ which can be viewed as 
a pluri-harmonic function in $\Sigma^{(l,m)}\times M$. Thus,
\[
\pi_2^*\omega + \p \bar \p (\phi^{(m,\epsilon)} - h^{(l,m)})  = \pi_2^*\omega + \p \bar \p \phi^{(m,\epsilon)},
\qquad \qquad \forall\;\; (s,t,\cdot) \in \Sigma^{(l,m)} \times M.
\] 
Thus, we have
\[
\left(\pi_2^*\omega + \p \bar \p \phi^{(l,m,\epsilon)}\right)^{n+1} =    (\pi_2^*\omega + \p \bar \p ( \phi^{(m,\epsilon)} - h^{(l,m)} ))^{n+1}
\]
with boundary data
\[
\begin{array}{lcl}
 \phi^{(l,m,\epsilon)}\mid_{\p \Sigma^{(l,m)} \times M} &= & \bar \rho \mid_{\p \Sigma^{(l,m)}\times M}
 \\ & \geq & \phi^{(m,\epsilon)} - h^{(l,m)}\mid_{\p \Sigma^{(l,m)} \times M} 
 \end{array}
\]
The last inequality holds because $\bar \rho= \phi^{(m,\epsilon)}$ when $t=0,m.\;$  In $\p \Sigma^{(l,m)}$
with $t\neq 0, m$, we have $$h^{(l,m)} = K > \phi^{(m)} - \bar \rho \geq \phi^{(m,\epsilon)}- \bar \rho. $$  
By the maximum principle for Monge-Ampere equation, we have
\[
\phi^{(m,\epsilon)} - h^{(l,m)} \leq \phi^{(l,m,\epsilon)}.
\]
In particular, we have
\begin{equation}
\phi^{(m)} +o(\epsilon) - h^{(l,m)} \leq \phi^{(l,m,\epsilon)} \leq \phi^{(m)}.
\end{equation}
Thus,  the sequence of K\"ahler potentials
$\phi^{(l, \epsilon)}$ converges to $\phi^{(m)}$ in  $\Sigma^{(\delta_l, m)} \times M\;$
in the $C^{0}$ norm.
\end{proof}

In fact, more is true.

\begin{lem} For any $\alpha \in (0,1), $ the sequence of K\"ahler potential
$\phi^{(l,m,\epsilon)}$ converges to $\phi$ in any sub-domain $\Sigma^{(\delta_l,m)} \times M\;$
in  the $C^{1,\alpha}$ norm.
\end{lem}

This is just a corollary of Theorem 5.2 and Lemma 5.11. We are now ready to prove Theorem 5.2.\\

 To obtain uniform $C^{1,1}$ bound of $\{\phi^{(l,m,\epsilon)}\}$ independent of $l, m \rightarrow \infty$
 and $\epsilon \rightarrow 0 ,\;$ we need to choose some appropriate background K\"ahler metric first. Let $h$ be the ambient K\"ahler metric which tame the initial geodesic ray $\rho: [0,\infty)\rightarrow \cH.\;$ Then,  the Dirichlet boundary value problem eq. \ref{eq:hcma0} can be re-written as
  a Dirichlet problem on
  $\Sigma^{(l,m)} \times M$ such that
\begin{eqnarray}
\det\left( h_{\alpha \b \beta} + {{\p^2 (\phi- \bar \rho)}\over {\p w^{\alpha}\p w^{\b \beta}}}\right)_{(n+1)\times (n+1)}
= \epsilon \cdot \det \left(h_{\alpha \b \beta}\right)_{(n+1)\times (n+1)}, \label{eq:hcma21}
\end{eqnarray}
with boundary condition
\begin{equation}
  \phi\mid_{\p \Sigma^{(l,m)} \times M} =  \bar \rho \mid_{\p \Sigma^{(l,m)}\times M}.
  \label{eq:hcmaT21}
\end{equation}

Lemma 5.11  implies a uniform $C^0$ bound on the solution $\phi^{(l,m,\epsilon)} - \phi^{(m)}\;$ where $\bar \rho$ is
a sub-solution of $\phi^{(l,m,\epsilon)}.\;$  To obtain a super-solution, we define a sequence
of (essentially) harmonic functions $\bar \psi^{(l,m)}$ by solving
\[
\tilde \triangle\; (\bar \psi^{(l,m)} -\bar \rho) + (n+1) = 0, \qquad \qquad \forall\;\;(s,t,\cdot)\in \Sigma^{(l,m)}\times M
\]   
with boundary data
\[
\b \psi^{(l,m)} \mid_{\p \Sigma^{(l,m)}\times M} = \bar \rho  \mid_{\p \Sigma^{(l,m)}\times M}. 
\]
Here $\tilde \triangle$ is the Laplacian operator of the K\"ahler metric $h.\;$
Since $h$ is a smooth metric with uniform bound on curvature and injectivity radius in
$\Sigma^{(l,m)}\times M$ and $\psi$ is a smooth function with uniform bounds (independent
of $l$), the standard elliptic PDE theory (interior and boundary estimate for harmonic
function)   implies that there is a uniform bound
on $\b \psi^{(l,m)}\;$ (up to two derivatives for instance).  Note that
\[
   \bar \rho \leq \phi^{(l,m,\epsilon)} \leq \bar \psi^{(l,m)}
\]
with equality holding on $\p \Sigma^{(l,m)}\times M.\;$  Consequently, we have the following
\begin{lem}  The first derivatives of $\phi^{(l,m,\epsilon)} -\bar \rho $ in $\p \Sigma^{(l,m)}\times M$ are uniformly bounded
(independent of $l$).
\end{lem}

 We want to solve equation (\ref{eq:hcma21}) and (5.11)  for any large $l, m \gg 1.\;$ Following \cite{chen991}, we have

\begin{lem}\label{cor:2ndderivativeestimate1} There exists a constant $ C $ which depend only on the ambient metric $ h$  ( but independent of $l,m$) and initial K\"ahler potentials $\varphi_0, \bar \rho (0) $ such that either
\[
e^{- (\rho -\b \rho)}
( n + 1+ \tilde{\triangle}(\phi^{(l,m,\epsilon)}-\bar \rho) ) \leq C
\]
or 
\[
e^{- (\rho -\b \rho)}
( n + 1+ \tilde{\triangle} (\phi^{(l,m,\epsilon)}-\bar \rho) ) \leq C \cdot \displaystyle \max_{\p \Sigma^{(l,m)} \times M} \; e^{-(\rho -\b \rho)}
( n + 1+ \tilde{\triangle} (\phi^{(l,m,\epsilon)}-\bar \rho) ). \]
Here we assume $\rho-\b \rho$ is positive since it is uniformly bounded anyway.
\end{lem}
  Using
the same boundary estimate as in  \cite{chen991}, we can obtain a relative $C^{1,1}$ estimate. Theorem 5.2 is
then proved. 

\begin{theo}  There exists a  uniform constant $C$ such that
\[
(n+1) + \triangle_{h} (\phi^{(l,m,\epsilon)}-\bar \rho) \leq  C.\;\] 

\end{theo}

\section{Future Problems}
Theorem 1.7 is  particularly interesting in the fano setting since we already know that  the existence
 of KE metric is equivalent to the properness of the K-energy functional. 
 This naturally leads to an old conjecture of G. Tian where Tian conjectured that the K energy functional
 is proper if and only if there exists a cscK metric in $[\omega].\;$   To better understand K\"ahler geometry from the point of view of this infinite dimensional space geometry, one may ask  the following conjecture/question.
   
\begin{conj} The existence of cscK metrics implies that the K-energy functional is proper in the sense that it bounds the  geodesic distance.
\end{conj}

Note that this version of properness is different from the one in Tian's original conjecture.
Tian's conjecture is perhaps more elusive.   There is a partial result in this direction which
 is a corollary of Theorem 1.1.

  \begin{cor} Let $\varphi_{0}$ be a cscK metric and let $\phi$ be any K\"ahler potential such
  that $\omega_{\phi}$ is uniformly elliptic.  If the geodesic distance of $\phi$ to any holomorphic line
  passing thorough $\varphi_{0}$ is greater than a fix constant $\delta > 0$, then the K energy of $\phi$ is bounded
  from below by a constant $C$ which depends only on $\delta$ and elliptic constant of $\omega_{\phi}.\;$
  \end{cor}
  
 In proof of  Theorem 1.1 and Corollary 6.2, we need to prove partially an old conjecture of this author.

\begin{conj} A global $C^{1,1}$ K-energy minimizer in any K\"ahler class must be smooth.
\end{conj}

This conjecture is proved in the Canonical K\"ahler class via weak K\"ahler Ricci flow.  \\

Another intriguing corollary of Theorem 1.1 is the following
  
\begin{cor} {\bf cscK metrics implies geodesic stable.}   If there is a cscK metric in the K\"ahler class, then $\yen$ invariant of every geodesic
ray must be strictly positive unless it is parallel to  a holomorphic line initiated from a cscK metric.
\end{cor}

Recently,  J. Stoppa \cite{stopp08} proves  a beautiful theorem in algebraic setting which fits into this discussion very well. Assuming there is no holomorphic vector fields, he proves that the existence of a cscK metric implies that the underlying polarization is K-stable. This is very reminiscent to a famous work of G. Tian \cite{tian97} where Tian proved that the existence of KE metric implies its polarization is K-stable with respect to special degenerations.  We also would like to note that T. Mabuchi announced that he can remove this assumption on automorphism group, again in algebraic case. 
  We will omit proofs to both corollaries since they basically  follow directly from the proof of Theorem 1.1.

 \begin{conj} For any destabilized test configuration such that the Calabi energy of nearby fibre approaches $0, $  if the total space is smooth or if the central fibre only admit  only mild singularities,  and if  the K-energy  of K\"ahler metrics of nearby fibres are uniformly bounded
 from below, then the corresponding polarization of nearby fibres has uniform lower bound on 
 K energy and furthermore it  is Semi-K-stable (whenever applicable).
 \end{conj}
 
 The first important case  to test is to show that Theorem 1.7 holds  when:  a) the Riem. curvature of the total space has a uniform lower bound ;
b)   the K-energy over the family of fibre metrics has a uniform lower bound; c) the Calabi energy
over the family of fibre metrics approaches to $0$.   In fact,
 the  author strongly suspects that semi-geodesic stability (perhaps also semi-K stability) will imply the existence of a uniform lower bound of the K-energy functional. \\
 
 While we are in this topic, the following question is very interesting.
 
 \begin{conj} In any K\"ahler manifold, if the K\"ahler Ricci flow (when applicable) or the Calabi flow
 initiated from one metric, converges to a cscK metric (by sequence geometrically), then the K-energy functional must have a lower bound in the original K\"ahler class. 
 \end{conj}
 
 Coming back to the variational nature of finding cscK metrics, the following question is
 interesting.
   
\begin{conj} For any {\bf critical} sequence of K\"ahler metrics,  can we replace with another {\bf critical} sequence where
the curvature is uniformly bounded away from a measure $0$ set? is this singular set has codimension 4? 
\end{conj}



\begin{thebibliography}{10}

\bibitem{Are-Tian03}
C. Arezzo and G. Tian.
\newblock Infinite geodesic rays in the space of KŠhler potentials.
\newblock  Ann. Sc. Norm. Super. Pisa Cl. Sci. (5) 2 (2003), no. 4, 617--630.

\bibitem{Bedford76}
E.D. Bedford and T.A. Taylor.
\newblock The {D}irichlet problem for the complex {M}onge-{A}mpere operator.
\newblock {\em Invent. Math.}, 37:1--44, 1975.

\bibitem{SongBen040}
B. Weinkove and J. Song.
\newblock On the convergence and singularities of the J-flow with applications to the Mabuchi energy
\newblock math.DG/0410418.

\bibitem{Ben04} 
Ben Weinkove.
\newblock  On the J-flow in higher dimensions and the lower boundedness of the Mabuchi energy
 \newblock math.DG/0309404.
 
\bibitem{CNS84}
 L.~Caffarelli, L.~Nirenberg and J.~Spruck.
\newblock The Dirichlet problem for nonlinear second-order elliptic equation
  {I}. Monge-Ampere equation.
\newblock {\em Comm. on pure and appl. math.}, XXXVII:369--402, 1984.

\bibitem{calabi82}
E.~Calabi.
\newblock Extremal {K}\"ahler metrics.
\newblock In {\em Seminar on Differential Geometry}, volume~16 of {\em 102},
  pages 259--290. Ann. of Math. Studies, University Press, 1982.

\bibitem{calabi85}
E.~Calabi.
\newblock Extremal {K}\"ahler metrics, {I}{I}.
\newblock In {\em Differential geometry and Complex analysis}, pages 96--114.
  Springer, 1985.

\bibitem{chen992}
E.~Calabi and X.~X. Chen.
\newblock Space of {K}\"ahler metrics and {C}alabi flow, 2002.
\newblock Journal of Differential Geometry.

\bibitem{chen991}
X.~X. Chen.
\newblock Space of {K}\"ahler metrics.
\newblock {\em Journal of Differential Geometry}, 56(2):189--234, 2000.

\bibitem{chen05}
X.~X. Chen.
\newblock Space of K\"ahler metrics (III)--the greatest lower bound of the Calabi energy.
\newblock to appear in Inventiones.


\bibitem{chen00}
X. ~X. Chen.
\newblock On the lower bound of the {M}abuchi energy and its application.
\newblock  {\em Internat. Math. Res. Notices} 2000, no. 12, 607--623. 

\bibitem{chensong081}
X.~X. Chen and S. Sun.
\newblock Private communication

\bibitem{chensong082}
X.~X. Chen and S. Sun.
\newblock Lower bound of the K-energy in algebraic manifolds (I).
\newblock In preparation.

\bibitem{chentang07}
X. X. Chen and Y. D. Tang.
\newblock test configuration and geodesic rays.
\newblock  arXiv:0707.4149 


\bibitem{chentian0} X. X. Chen, G. Tian.
\newblock {\it Ricci flow on K\"ahler-Einstein surfaces},
\newblock Invent. math. {\bf 147} (2002), 487-544

\bibitem{chentian005}
X.~X. Chen and G.~Tian.
\newblock Foliation by holomorphic discs and its application in K\"ahler
  geometry, 2003.
\newblock to appear in Publication in Mathematics.

\bibitem{chenhe05}
X.~X. Chen and W. Y. He.
\newblock On the Calabi flow,
\newblock Amer. J. Math. 130 (2008), no. 2, 539--570.


\bibitem{Semmes93}
R.R. Coifman and S.~Semmes.
\newblock Interpolation of {B}anach spaces, {P}erron process, and {Y}ang
  {M}ills.
\newblock {\em Amer. J. Math.}, 115(2):243--278, 1993.

\bibitem{Dona96}
S.~K. Donaldson.
\newblock Symmetric spaces, {K}{\"a}hler geometry and {H}amiltonian dynamics.
\newblock {\em Amer. Math. Soc. Transl. Ser. 2, 196}, pages 13--33, 1999.
\newblock Northern California Symplectic Geometry Seminar.

\bibitem{Dona01}
S.~K. Donaldson.
\newblock Holomorphic discs and the complex {M}onge-{A}mp\`ere equation, 2001.
\newblock to appear in {J}ournal of {S}ympletic {G}eometry.

\bibitem{Dona05}
S. K. Donaldson.
\newblock Scalar curvature and projective embedding. {II}. 
\newblock {\em Q. J. Math.} 56 (2005), no. 3, 345--355. 
 	
\bibitem{Dona051}
S. K. Donaldson.
\newblock  Interior estimates for solutions of Abreu's equation.
\newblock  {\em Collect. Math.}  56 (2005), no. 2, 103--142.

\bibitem{Dona052}
S. K. Donaldson.
\newblock Lower bounds on the Calabi functional.
\newblock  J. Differential Geom. 70 (2005), no. 3, 453--472. 
 
\bibitem{Dona081}
S. K. Donaldson.
\newblock  A note  on the $\alpha$ invariant  of Mukai-Umermura 3-fold
\newblock  arXiv: 0711.3357v1. 
 

\bibitem{Futaki83}
A. Futaki.
\newblock An obstruction to the existence of Einstein K\"ahler 
metrics.
\newblock  {\em Invent. Math.} 73 (1983), no. 3, 437--443.

\bibitem{FutakiMa95}
A. Futaki and T. Mabuchi.
\newblock  Bilinear forms and extremal K\"ahler vector fields associated with K\"ahler classes. 
\newblock {\em Math. Ann.} 301 (1995), no. 2, 199--210.

\bibitem{GuanB98}
B. Guan.
\newblock The {D}irichlet problem for complex {M}onge-{A}mpere equations and
  regularity of the plui-complex green function.
\newblock {\em Comm. {A}na. {G}eom.}, 6(4):687--703, 1998.


\bibitem{Hwang951}
A.~D. Hwang.
\newblock On the {C}alabi energy of {E}xtremal {K}\"ahler metrics.
\newblock {\em International Journal of Mathematics}, 6(6):825--830, 1995.


\bibitem{CKNS85}
J.T. Kohn, L.~Nirenberg,  L.~Caffarelli and J.~Spruck.
\newblock The Dirichlet problem for nonlinear second-order elliptic equation
  {II}. complex monge-ampere equation.
\newblock {\em Comm. on pure and appl. math.}, 38:209--252, 1985.

\bibitem{Ma87}
T.~Mabuchi.
\newblock Some symplectic geometry on compact {K}\"ahler manifolds {I}.
\newblock {\em Osaka, {J}. {M}ath.}, 24:227--252, 1987.

\bibitem{Ma05}
T. ~Mabuchi. 
\newblock An energy-theoretic approach to the Hitchin-Kobayashi correspondence for manifolds. {I}. 
\newblock  {\em Invent. Math.} 159 (2005), no. 2, 225--243.

\bibitem{paultian04}
 S. Paul and G. Tian. 
 \newblock Analysis of geometric stability. 
 \newblock {\em Int. Math. Res. Not.} 2004, no. 48, 2555--2591. 
 
 \bibitem{paul041}
S. Paul.
\newblock Geometric analysis of Chow Mumford stability. 
\newblock {\em Adv. Math.} 182 (2004), no. 2, 333--355.


\bibitem{Phong_Sturm05}
D.H. Phong and J. Sturm.
\newblock The Monge-Amp\`ere operator and geodesics in the space of K\"ahler potentials.
\newblock math.DG/0504157.

\bibitem{ruan90}
W.D.~Ruan.
\newblock  On the convergence and collapsing of {K}\"ahler metrics.
\newblock   {\em J. Diff. Geom.}, 52:1-40, 1999.

\bibitem{Semmes92}
S.~Semmes.
\newblock Complex {M}onge-{A}mp\`ere equations and sympletic manifolds.
\newblock {\em Amer. J. Math.}, 114:495--550, 1992.


\bibitem{song08}
S. Sun.
\newblock Note on geodesic rays tamed by simple test configurations.
\newblock preprint, arXiv:0806.2697.

\bibitem{stopp08}
J. Stoppa.
\newblock K-stability of constant scalar curvature KŠhler manifolds.
\newblock preprint, arXiv:0803.4095.
\bibitem{tian87}
G. Tian.
\newblock On K\"ahler-Einstein metrics on certain K\"ahler manifolds with $c_1(M)>0.$
\newblock {\em Invent. Math.} 89 (1987),  225-245.

\bibitem{tian97}
G. Tian.
\newblock K\"ahler-Einstein metrics with positive scalar curvature. 
\newblock {\em Invent. Math.} 130 (1997), no. 1, 1--37. 

\bibitem{Yau78}
S.~T. Yau.
\newblock On the {R}icci curvature of a compact {K}\"ahler manifold and the
  complex {M}onge-{A}mpere equation, ${I}^*$.
\newblock {\em Comm. Pure Appl. Math.,}, 31:339--441, 1978.

\end{thebibliography}
\end{document}